\documentclass[12pt]{amsart}
\usepackage{amsmath,amssymb, mathtools} 
\usepackage{amsthm} 
\usepackage{dsfont} 
\usepackage{tensor} 

\usepackage[pdftex]{graphicx} 
\usepackage{color}
\usepackage[colorlinks,citecolor=blue,linkcolor=blue,urlcolor=blue]{hyperref} 
\usepackage{cleveref}

\usepackage[utf8]{inputenc} 
\usepackage[english]{babel} 

\usepackage{geometry} 
\usepackage{tikz} 
\usepackage{cite} 
\usepackage{enumitem} 

\usepackage{mathrsfs}

\newtheorem{theorem}{Theorem}

\newtheorem{lemma}[theorem]{Lemma}

\newtheorem{proposition}[theorem]{Proposition}
\newtheorem{example}[theorem]{Example}
\newtheorem{Remark}[theorem]{Remark}
\newtheorem*{procedure}{Procedure}

\numberwithin{theorem}{section} 
\numberwithin{equation}{section} 
\numberwithin{figure}{section} 

\newcommand*{\mf}{\mathfrak}
\newcommand*{\mc}{\mathcal}
\newcommand*{\bb}{\mathbb}
\newcommand*{\R}{\mathbb{R}}
\newcommand*{\C}{\mathbb{C}}

\newcommand*{\g}{\mathfrak{g}}

\DeclareMathOperator{\tr}{tr}
\DeclareMathOperator{\ad}{ad}
\DeclareMathOperator{\SL}{SL}
\DeclareMathOperator{\Hom}{Hom}

\def\bbone{{\mathchoice {\rm 1\mskip-4mu l} {\rm 1\mskip-4mu l}
{\rm 1\mskip-4.5mu l} {\rm 1\mskip-5mu l}}}

\title{Structure constants for simple Lie algebras from principal $\mf{sl}_2$-triple}

\author{Abdelmalek Abdesselam}
\author{Alexander Thomas}
\address{Abdelmalek Abdesselam, Department of Mathematics,
P. O. Box 400137,
University of Virginia,
Charlottesville, VA 22904-4137, USA}
\address{Alexander Thomas, Mathematikon, Universit\"at Heidelberg, Berliner Str. 41-49, 69120 Heidelberg, Germany}
\email{malek@virginia.edu}
\email{athomas@mathi.uni-heidelberg.de}
\date{}
	
\begin{document}

\maketitle

\begin{abstract}
For a simple complex Lie algebra $\g$, fixing a principal $\mf{sl}_2$-triple and highest weight vectors induces a basis of $\g$ as vector space. For $\mathfrak{sl}_n(\C)$, we describe how to compute the Lie bracket in this basis using transvectants. This generalizes a well-known rule for $\mf{sl}_2$ using Poisson brackets and degree 2 monomials in two variables. Our proof method uses a graphical calculus for classical invariant theory. Other Lie algebra types are discussed.
\end{abstract}

\section{Introduction and results}

\subsection{The rank 1 case}

Consider the Lie algebra $\mathfrak{sl}_2$ in its standard representation as traceless $2\times 2$-matrices (all Lie algebras are over $\C$). Fix the standard basis $(E,H,F)$ where
$$E=\begin{pmatrix} 0 & 1\\ 0 & 0\end{pmatrix}, \;\; H=\begin{pmatrix} 1 & 0\\ 0 & -1\end{pmatrix}, \;\; F=\begin{pmatrix} 0 & 0\\ 1 & 0\end{pmatrix}.$$

The Lie bracket relations read $[H,E]=2E, [H,F]=-2F$ and $[E,F]=H$. We can recover these relations in the following way: Using the identification \begin{equation}\label{sl2-identity}
(E,H,F)=\left(\tfrac{x^2}{2},-xy,-\tfrac{y^2}{2}\right),
\end{equation}
the Lie bracket simply becomes the Poisson bracket $\{.,.\} = \frac{\partial}{\partial x}\wedge \frac{\partial}{\partial y}$. For example we have
$\{-xy,-y^2/2\}=y^2$ which mimics $[H,F]=-2F$.

Why does this method work? First of all, $\mf{sl}_2$ acts on $\C^2$, so on its function space $\C[x,y]$ preserving the degree. On polynomials of fixed degree $d$ we get the (unique) irreducible representation of dimension $d+1$. The adjoint action being of dimension 3, it is given by the action on polynomials of degree 2.

The link to the Poisson bracket comes from the fact that $\mf{sl}_2\cong \mf{sp}_2$ and $\mathrm{Sp}_2(\C)$ acts on $\C^2$ equipped with the complex symplectic structure $dx\wedge dy$ in a Hamiltonian way. This means that the vector field associated to the infinitesimal action of $x\in\mf{sp}_2$ is the symplectic gradient of a function $H_x$ such that $H_{[x,y]}=\{H_x,H_y\}$. The functions $H_x$ for $x\in\{E,H,F\}$ are the monomials from Equation \eqref{sl2-identity}.

The question arises whether a similar rule holds for $\mathfrak{sl}_n$ and other simple Lie algebras.
The action of $\mathrm{Sp}_{2n}(\C)$ on $\C^2$ is Hamiltonian, so we can generalize the correspondence \eqref{sl2-identity} via Poisson brackets to $\mathfrak{sp}_{2n}$. The drawback is that the polynomials associated to the generators of $\mathfrak{sp}_{2n}$ are not monomials anymore. 

Another direction for generalization, and this is the route we take here, is to use a higher order version of the Poisson bracket: the \emph{transvectants} 
of classical invariant theory (see, e.g.,~\cite{gordan}). They are $\mathrm{SL}_2(\R)$-invariant bilinear differential operators on the circle. Transvectants appear in the Clebsch--Gordan decomposition and in the Moyal product, the deformation of the Poisson bracket on $\R^{2n}$. We give more details in Section \ref{Sec:graphical-calculus} around Equation \eqref{Eq:transvectant-graphical}, and also recommend \cite[Chapter 3.1]{ovsienko2004projective} for an introduction to transvectants from the point of view of Poisson geometry and the theory of integrable systems.

\subsection{Generalization}

We describe our procedure in the case of Lie algebras of type $A$, i.e. for $\mathfrak{sl}_n(\C)$. For other types, see Section \ref{Sec:other-types}.

Fix a principal $\mathfrak{sl}_2$-triple $(E,H,F)$ of $\mathfrak{sl}_n(\C)$. Thus, we can consider $\mathfrak{sl}_n$ as $\mathfrak{sl}_2$-module and we can decompose it into irreducible representations. In the standard representation of $\mathfrak{sl}_n$ on $\C^n$, a natural choice of highest weight vectors is given by $(E^k)_{1\leq k\leq n-1}$. Then we get a basis of $\mathfrak{sl}_n$ as vector space by applying $\ad_F$ successively to the highest weight vectors. Denote this basis by $G_{i,j}$, where $G_{i,j}=\ad_F^{i-j}(E^i)$ with $i\in \{1,...,n-1\}$ and $j\in \{-i,-i+1,...,i-1,i\}$.

To $G_{i,j}$ we associate the monomial 
\begin{equation}\label{Eq:g-monomial}
m_{i,j} = y^{i-j}\frac{x^{i+j}}{(i+j)!}.
\end{equation}
 
\begin{procedure}
Computing the Lie bracket between two elements in the basis $G_{i,j}$ goes as follows:
\begin{enumerate}
    \item Translate the two elements of $\mathfrak{sl}_n$ into two polynomials using \eqref{Eq:g-monomial},
    \item Compute the transvectants of the two polynomials and sum them with specific weights,
    \item Translate the result to $\mathfrak{sl}_n$ using \eqref{Eq:g-monomial} backwards.
\end{enumerate}
\end{procedure}

In the second step, the $m$-th transvectant between the monomials $m_{k,i}$ and $m_{\ell,j}$ has to be weighted by a constant $\omega_{k,\ell}^{(m)}$ (independent of $i$ and $j$, and where $k,\ell\in\{1,...,n-1\}$ and $m\in\{1,...,2\,\mathrm{min}(k,\ell)-1\}$).
We call them the \emph{structure constants of the Lie bracket} in $\mathfrak{sl}_n$. 

To compute the transvectant of two monomials $x^ay^b$ and $x^cy^d$, define $\Omega$ to be the following map on $\C[x,y]\otimes \C[x,y]$ to itself:
$$\Omega(x^ay^b\otimes x^cy^d) = ax^{a-1}y^b\otimes dx^cy^{d-1}-bx^ay^{b-1}\otimes cx^{c-1}y^d.$$
This map $\Omega$ ``applies the Poisson bracket''. To compute the $m$-th transvectant, iterate $m$ times $\Omega$ and finally compose with the multiplication map $\C[x,y]\otimes\C[x,y]\to \C[x,y]$. This process is the special case of \emph{Cayley's $\Omega$-process} for binary forms. To be explicit, for $F,G\in \C[x,y]$ homogeneous of degree $k$ and $\ell$ respectively and for $0\leq m\leq \mathrm{min}(k,\ell)$, the $m$-th transvectant $(F,G)_m$ is a homogeneous polynomial of degree $k+\ell-2m$ given by:
\begin{equation}
(F,G)_m=
\frac{(k-m)!\ (\ell-m)!}{k!\ \ell!}\ 
\sum_{j=0}^m
(-1)^j
\binom{m}{j}
\frac{\partial^m F}{\partial x^{m-j}\partial y^{j}}
\frac{\partial^m G}{\partial x^{j}\partial y^{m-j}}\ .
\end{equation}

From this expression, we see the symmetry $(F,G)_m=(-1)^m(G,F)_m$, which implies $\omega_{k,\ell}^{(m)}=(-1)^{m+1}\omega_{\ell,k}^{(m)}$. In particular, the structure constants vanish for $m$ even, a result which we recover in Theorem \ref{Thm:2} as well as in Proposition \ref{Prop:vanishing}. This is why we often restrict to $m$ odd, in which case the structure constant are symmetric in exchanging $k$ and $\ell$.

\begin{example}
For $\mathfrak{sl}_3$, we can compute $[E^2,F^2]$ with our procedure:
$$\left(\frac{x^4}{4!},\frac{y^4}{4!}\right)\mapsto \omega_{2,2}^{(1)}\frac{x^3y^3}{3!3!}+\omega_{2,2}^{(3)}xy.$$
For $\mathfrak{sl}_3$, there are no monomials of degree 6 (see Figure \ref{Fig:monomial-asso}), so $\omega_{2,2}^{(1)}=0$, and we have $\omega_{2,2}^{(3)}=-2$. Since $xy$ corresponds to $-H$, we get $[E^2,F^2]=2H$.

For $\mathfrak{sl}_4$, we have $\omega_{2,2}^{(1)}=2/5$ and $\omega_{2,2}^{(3)}=-24/5$ (see Appendix \ref{Sec:appendix}). Further, $\frac{x^3y^3}{3!3!}$ corresponds to $3!\,G_{3,0}$. Hence $[E^2,F^2]=12/5\, G_{3,0}-24/5\,H$.
\end{example}

\subsection{Results}

Our first theorem states that the procedure above works. Let $V_i=\mathrm{Span}(G_{i,j}\mid -i\leq j\leq i)$ and $\mathcal{H}_{2i}=\C_{2i}[x,y]$, the space of homogeneous polynomials of degree $2i$. Denote by $\Phi_i:V_i\to \mathcal{H}_{2i}$ the vector space isomorphism sending $G_{i,j}$ to the monomial $m_{i,j}$ from \eqref{Eq:g-monomial}.

\begin{theorem}\label{Thm:1}
    There exist constants $\omega_{k,\ell}^{(m)}(\mathfrak{sl}_n)$ with $1\leq \ell\leq k\leq n-1=\mathrm{rk}(\mathfrak{sl}_n)$ and $1\leq m\leq 2\ell-1$, such that the above procedure computes the Lie bracket in $\mathfrak{sl}_n$. That is, for all $1\leq \ell\leq k\leq n-1$ and $-k\leq i\leq k$ and $-\ell\leq j\leq \ell$ we have:
$$[G_{k,i},G_{\ell,j}]=\sum_{m=1}^{2\ell-1} \omega_{k,\ell}^{(m)}(\mathfrak{sl}_n)\,\Phi_{k+\ell-m}^{-1}((m_{k,i},m_{\ell,j})_m).$$
\end{theorem}
We will often drop the dependence of $\mathfrak{sl}_n$ in the sequel since we work for fixed $n$.  The values of the structure constants for $\mathfrak{sl}_n$ for $n\leq 6$ are given in the Appendix \ref{Sec:appendix}. Note that $\omega_{k,1}^{(m)}=\delta_{1,m}$, so we can consider $\ell>1$.

Our main result computes the structure constants explicitly for $\mathfrak{sl}_n$. Put
$$\mathcal{Q}_{k,\ell}^{(m)}:=2\times(-1)^{k+\ell+n-1}\times(2k+2\ell-2m+1)\times
\frac{k!^2\ \ell!^2\ (n-k-\ell+m-1)!}{m!\ (n+k+\ell-m)!}\ ,$$
and
$$\mathcal{R}_{k,\ell}^{(m)}:=
\sum_{q\in\mathbb{Z}}(-1)^q
\binom{q\! +\! 1}{2k\! +\! 2\ell \! -\! m\! +\! 1}\!
\binom{m}{q\! -\! k\! -\! \ell\! +\! m\! -\! n\! +\! 1}\!
\binom{2k\! -\! m}{q\! -\! \ell\! -\! n\! + \! 1}\!
\binom{2\ell\! -\! m}{q\! -\! k\! -\! n\! +\! 1}\ ,$$
with the convention that binomials are defined as zero if lying outside of Pascal's triangle. It is obvious from the above formulas that $\mathcal{Q}_{\ell,k}^{(m)}=\mathcal{Q}_{k,\ell}^{(m)}$ and $\mathcal{R}_{\ell,k}^{(m)}=\mathcal{R}_{k,\ell}^{(m)}$.

\begin{theorem}\label{Thm:2}
    For integers $n, k, \ell, m$ with $n\geq 2$, $1\leq \ell\leq k\leq n-1$ and $1\leq m\leq 2\ell-1$, we have
    $$\omega_{k,\ell}^{(m)}(\mathfrak{sl}_{n}) = 0,$$ if $m$ is even, and
    $$\omega_{k,\ell}^{(m)}(\mathfrak{sl}_{n}) = \mathcal{Q}_{k,\ell}^{(m)}\times \mathcal{R}_{k,\ell}^{(m)}, $$ if $m$ is odd.
\end{theorem}

In two special cases, the expression simplifies a lot:
\begin{itemize}
    \item For $m=1$, we simply get $\omega_{k,\ell}^{(1)}=\frac{1}{2}\frac{(2k)!(2\ell)!}{(2k+2\ell-2)!}$ which is independent of $n$.
    \item For $m=2\ell-1$, the maximal possible value of $m$, we get $$\omega_{k,\ell}^{(2\ell-1)} = (-1)^{\ell+1}\ell\times\frac{(k!)^2}{(k-\ell+1)!^2}\times\frac{\binom{n+k}{2k+1}}{\binom{n+k-\ell+1}{2k-2\ell+3}}.$$
\end{itemize}

The explicit expression gives some symmetries of the structure constants: as already mentioned, for $m$ odd there is a symmetry by exchanging $k$ and $\ell$: $\omega_{k,\ell}^{(m)} = \omega_{\ell,k}^{(m)}$. A more subtle symmetry is an explicit proportionality factor between $\omega_{k,\ell}^{(m)}$ and $\omega_{k,k+\ell-m}^{(2k-m)}$, see Proposition \ref{Prop:hidden-symmetry} below.

Formulas similar to the one in Theorem \ref{Thm:2} are common in the theory of angular momentum in quantum mechanics initiated by Wigner \cite{Wigner} and surveyed in the book \cite{biedenharn1984angular}. Our approach to prove Theorem \ref{Thm:2} is to use the graphical calculus developed by the first author in \cite{AbdesselamJKTR}, unifying the classical invariant theory of binary forms and the theory of angular momentum. An expression of the structure constants $\omega_{k,\ell}^{(m)}$ in terms of Wigner $6j$-symbols is given in Equation \eqref{Eq:6j-expr}.

Another source for formulas similar to the one in Theorem \ref{Thm:2} is the theory of $W$-algebras. In particular, the $W_\infty$-algebra defined in \cite{pope1990w} has similar structure constants. There might be a link between our structure constants for $\mathfrak{sl}_n$ in some suitable limit $n\to \infty$ and the structure constants for $W_\infty$.

\subsection{Application: character varieties}

The initial motivation for this project comes from the study of character varieties of surface groups and in particular the Hitchin component \cite{HCS}.

Consider a Riemann surface $S$ and a complex simple Lie group $G$. The character variety $\mathrm{Rep}(\pi_1S, G)$ is the space of all isomorphism classes of completely reducible representations of $\pi_1S$ into $G$:
$\mathrm{Rep}(\pi_1S, G)=\mathrm{Hom}(\pi_1S,G)/G,$
where $G$ acts by conjugation.
These character varieties appear in many contexts, for instance geometric structures, flat connections (via the Riemann--Hilbert correspondence) or moduli spaces of holomorphic bundles (via the non-abelian Hodge correspondence). Probably the most important example is \emph{Teichmüller space}, the moduli space of complex structures on a surface up to isotopy, which is a connected component of $\mathrm{Rep}(\pi_1S, \mathrm{PSL}_2(\R))$.

Hitchin's approach to Teichmüller space in \cite{hitchin1992lie}, which leads to a generalization to higher rank, is to consider a special holomorphic bundle $V$ together with a holomorphic $(1,0)$-form $\Phi$, the so-called \emph{Higgs field}, which lies in the Kostant slice. From this data, Hitchin constructs a flat connection with monodromy in $\mathrm{PSL}_n(\R)$.

For a principal $\mathfrak{sl}_2$-triple $(E,H,F)$ in a simple Lie algebra $\g$, the \emph{Kostant slice} is the affine space $F+Z(E)$, where $Z(E)=\{x\in \g\mid [x,E]=0\}$ denotes the centralizer.
Kostant's slice theorem says that almost every element\footnote{These elements are called \emph{regular}.} in $\g$ can be conjugated in a unique way to an element of this slice \cite{kostant59}. In the standard representation of $\g=\mathfrak{sl}_n$, $Z(E)$ is generated as a vector space by $(E^k)_{1\leq k\leq n-1}$.

A more general framework for studying character varieties is the following: instead of working with a holomorphic vector bundle, we consider simply a complex vector bundle $V$ of degree 0 (to allow flat connections). On $V$ we consider a field $\Phi\in \Omega^1(S,\g)$, i.e. a $\g$-valued 1-form. Using the Hodge decomposition we can write $\Phi=\Phi^{1,0}+\Phi^{0,1}$. We impose $\Phi\wedge\Phi=0$ and we prescribe the conjugacy class of $\Phi^{0,1}$.
In Hitchin's setting, $\Phi^{0,1}=0$ so its conjugacy class is zero and $\Phi\wedge\Phi=0$ is automatically true.

Another interesting case is when the conjugacy class of $\Phi^{0,1}$ is the one of a principal nilpotent element. This condition turns out not to depend on the Riemann surface structure on $S$.
Given a principal $\mathfrak{sl}_2$-triple $(E,H,F)$ of $\mathfrak{sl}_n$, we can locally put $\Phi^{0,1}=Fd\bar{z}$, so $\Phi^{1,0}\in Z(F)dz$, where $(z,\bar{z})$ is a local coordinate system on $S$. A variation of this data changes the conjugacy class of $\Phi^{0,1}$. Up to conjugation, the variation of $\Phi^{0,1}$ lies in the Kostant slice $F+Z(E)$. 

Hence we are in a ``doubled'' setting of the Kostant slice: $\Phi^{0,1}\in F+Z(E)$ and $\Phi^{1,0}\in Z(\Phi_2)$. There are $2\,\mathrm{rk}(\g)$ degrees of freedom here: $\mathrm{rk}(\g)$ degrees of freedom in the Kostant slice for $\Phi^{0,1}$, and the same for the choice of $\Phi^{1,0}$.

In the attempt to construct a flat connection out of the data $(V,\Phi)$ there is a set of constraints on the $2\,\mathrm{rk}(\g)$ degrees of freedom, see \cite[Section 7]{HCS}. To compute these constraints the formulas here are useful since both highest weight vectors $E^k$ and lowest weight vectors $F^\ell$ appear.

\subsection{Structure of the paper}

In Section \ref{Sec:preliminaries} we expose some preliminary notions from Lie theory, and introduce the natural basis. Then in Section \ref{Sec:structure-constants} we discuss the structure constants. In Section \ref{Sec:trace-method} we study a first method using the trace computing some structure constants. The core part of the paper is Section \ref{Sec:graphical-calculus} which introduces the graphical calculus, leading finally to the complete proof of Theorem \ref{Thm:2}. In the final Section \ref{Sec:other-types} Lie algebras different from $\mathfrak{sl}_n$ are discussed.

\medskip
\textbf{\textit{Acknowledgements.}} We are grateful to the platform MathOverflow which brought us together and catalysed this collaboration. We warmly thank the referee for numerous useful comments. We acknowledge support from the University of Heidelberg and the University of Virginia. The second author was supported by the European Research Council under ERC-Advanced Grant 101018839 and by the Deutsche Forschungsgemeinschaft (DFG, German Research Foundation) - Project-ID 281071066 - TRR 191.

\section{Preliminaries}\label{Sec:preliminaries}


\subsection{Natural basis from principal $\mathfrak{sl}_2$-triple}\label{Sec:natural-basis}

Consider a simple complex Lie algebra $\g$. An \emph{$\mathfrak{sl}_2$-triple} in $\g$ is the image of an injective homomorphism of $\mathfrak{sl}_2$ into $\g$. It is called \emph{principal} if the image of any non-zero nilpotent element of $\mathfrak{sl}_2$ is principal nilpotent in $\g$, i.e. is a nilpotent element with minimal centralizer (of dimension equal to the rank of $\g$). By a theorem of Kostant \cite{kostant59}, we know that there is a unique principal $\mathfrak{sl}_2$-triple in $\g$ up to conjugation. Fix $(E,H,F)$ such a triple. For $\mathfrak{sl}_n$, one possible choice is given by
\begin{equation}\label{Eq:e-h-f}
    E=\textstyle\sum_{j=1}^{n-1} E_{j,j+1}\, , H=\sum_{j=1}^n (n-2j+1)E_{j,j} \;\;\text{ and } \;\; F=\sum_{j=1}^{n-1}r_jE_{j+1,j}
\end{equation}
where $E_{i,j}$ denotes the standard basis for matrices and $r_j=j(n-j)$.

The principal triple induces two decompositions of $\g$. First by weights of $\ad_H$:
$$\g \cong \bigoplus_{k\in\mathbb{Z}} \g_k \;\; \text{ where } \;\; \g_k=\{g\in \g\mid [H,g]=k g\}.$$
Second by the action with the bracket, $\g$ becomes an $\mathfrak{sl}_2$-module which can be decomposed into irreducible representations:
$$\g \cong \bigoplus_{i\in\mathbb{Z}_{>0}} n_iV_i$$
where $V_i$ is the irreducible representation of $\mathfrak{sl}_2$ of dimension $2i+1$ and $n_i\in \mathbb{N}$ the multiplicities.

\begin{figure}[h!]
\centering
\includegraphics[height=4cm]{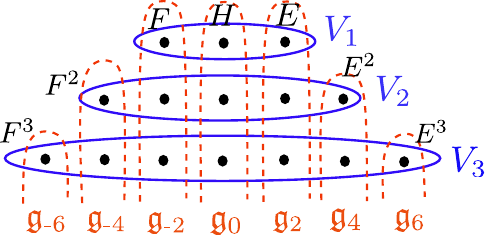}

\caption{Line decomposition of $\mathfrak{sl}_4$ by a principal triple}\label{Fig:g-decompo}
\end{figure}

From now on, we consider $\g=\mathfrak{sl}_n$. Then we know that $n_i\in \{0,1\}$. Using both decompositions, we get
\begin{equation}\label{line-decompo}
\mathfrak{sl}_n\cong \bigoplus \g_k\cap V_i,
\end{equation}
which is a line decomposition. See Figure \ref{Fig:g-decompo} for an illustration.

All finite-dimensional irreducible representations of $\mathfrak{sl}_2$ are highest weight representations. This means that for a given irreducible representation $V$, there is a vector $v\in V\backslash\{0\}$ with $E.v=0$, called highest weight vector. Then acting successively with $F$ generates all of $V$. 

In our setting, the highest weight vector of $V_i$ is given by $E^i$. Hence a basis adapted to the line decomposition \eqref{line-decompo} is given by
\begin{equation}\label{def-G}
G_{i,j}=\ad_F^{i-j}(E^i) \in V_i\cap \g_{2j}
\end{equation}
where $i\in \{1,...,n-1\}$ and $j\in \{-i,-i+1,...,i-1,i\}$.

There are two lowest weight vectors in $V_i$: $F^i$ and $G_{i,-i}=\ad_F^{2i}(E^i)$. Since there is a unique lowest weight vector up to scale, both differ by some constant given by the following proposition:
\begin{proposition}\label{prop:lowest-weight-link}
We have $G_{i,-i}=(-1)^i (2i)! \, F^i$.
\end{proposition}
The identity is independent of the choice of principal $\mathfrak{sl}_2$-triple, since all of them are conjugated. Then the proof reduces to a direct computation using the representatives \eqref{Eq:e-h-f}.

Similar to the case of $\mathfrak{sl}_2$, we can now associate monomials to the $G_{i,j}$. Then, we want to recover the Lie bracket as a manipulation of the polynomials. We fix the following association:
\begin{equation}
G_{i,j} \cong y^{i-j}\frac{x^{i+j}}{(i+j)!}.
\end{equation}

\begin{figure}[h!]
\centering
\includegraphics[height=4.5cm]{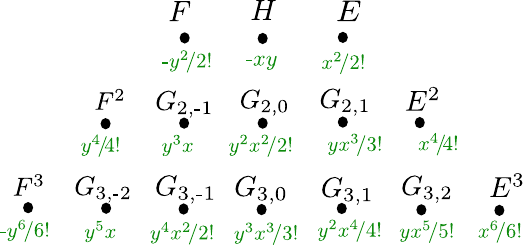}

\caption{Associated monomials to the line decomposition. }\label{Fig:monomial-asso}
\end{figure}

Note that the $\mathfrak{sl}_2$-triple $(E,H,F)$ gets associated to $(x^2/2,-xy,-y^2/2)$ which is slightly different, but equivalent to Equation \eqref{sl2-identity}. This is because $G_{1,0}=[F,E]=-H$ and $G_{1,-1}=\ad_F^2(E)=-2F$. Figure \ref{Fig:monomial-asso} shows the monomials for $G_{i,j}$ with $i\leq 3$. Note that the lowest weight vector $F^k$ is associated to $(-1)^ky^{2k}/(2k)!$ by Proposition \ref{prop:lowest-weight-link}.

\subsection{The Lie bracket}\label{Sec:Lie-bracket}

Recall that the Lie bracket $[.,.]:\g\times \g\to\g$ is bilinear, antisymmetric and satisfies the Jacobi identity. This implies that $[.,.]\in\Hom_{\mathfrak{sl}_2}(\Lambda^2\g,\g)$ where $\mathfrak{sl}_2$ acts on $\g$ via the principal triple which induces an action on $\Lambda^2\g$.
Restricted to $V_k\times V_\ell$, we get 
$$[.,.]\in \Hom_{\mathfrak{sl}_2}(V_k\otimes V_\ell, \oplus_i V_i) = \oplus_i \Hom_{\mathfrak{sl}_2}(V_k\otimes V_\ell, V_i).$$
We can decompose $V_k\otimes V_\ell$ into irreducible representations of $\mathfrak{sl}_2$ using the Clebsch-Gordan theorem. This gives for $k\geq \ell$:
$$V_k\otimes V_\ell \cong \oplus_{m=0}^{2\ell} V_{k+\ell-m}.$$
Finally, by Schur's lemma, we know that
$$\Hom_{\mathfrak{sl}_2}(V_m,V_i) \cong \delta_{i,m}\C.$$
Now, the projection of the Lie bracket $V_k\otimes V_{\ell}\to V_{k+\ell-m}$ given by $(v_k\otimes v_\ell)\mapsto \pi_{k+\ell-m}([v_k,v_\ell])$, where $\pi_{k}:\mathfrak{sl}_n\to V_k$ denotes the projection, is $\mathfrak{sl}_2$-equivariant. The same holds for the projection $V_k\otimes V_\ell \to V_{k+\ell-m}$ using the Clebsch--Gordan decomposition, which is given by the $m$-th transvectant. Since $\mathrm{dim}_{\mathfrak{sl}_2}(V_k\otimes V_\ell,V_{k+\ell-m})\leq 1$, the two maps are proportional. We denote by $\omega_{k,\ell}^{(m)}$ the proportionality factors, which represent the \emph{structure constants} of the Lie bracket in the basis $G_{i,j}$.

Since $V_k$ is stable under the adjoint action of $V_1$ (which is the principal triple), and since $G_{k,k-1}=[F,E^k]$, we get $\omega_{k,1}^{(m)}=\delta_{1,m}$. So we concentrate on $\ell>1$ in the sequel.
The bracket between $G_{i,j}$ and $G_{k,\ell}$ is a sum of elements in $\g_{2j+2\ell}$ which represents a column in Figure \ref{Fig:g-decompo}.

\section{The structure constants}\label{Sec:structure-constants}

The proof of Theorem \ref{Thm:1}, that our procedure for computing the Lie bracket works, is now easy.

\begin{proof}[Proof of Theorem \ref{Thm:1}]
    Consider two elements $x,y$ of $\mathfrak{sl}_n$ for which we want to compute $[x,y]$. Since our procedure respects the bilinearity of the two entries, we can suppose $x=G_{i,j}\in V_i$ and $y=G_{k,\ell}\in V_k$ for some integers $i,j,k,\ell$.

    In the previous subsection \ref{Sec:Lie-bracket}, we have seen that the Lie bracket restricted to $V_i\otimes V_k$ is an element of $\Hom_{\mathfrak{sl}_2}(V_i\otimes V_k, \mathfrak{sl}_n)$. The projection onto the factor $V_j$ of $\mathfrak{sl}_n$ is then given by some multiple of the projection from $V_i\otimes V_k$ onto $V_j$ which is given by the transvectant of order $i+k-j$. This multiple only depends on the decomposition of $\mathfrak{sl}_n$ into irreducible $\mathfrak{sl}_2$-modules. 

    Since $\mathfrak{sl}_n\cong \oplus_{j=1}^{n-1} V_j$, we get the result of $[x,y]$ by summing over all transvectants, weighted by the structure constants.
\end{proof}

We can now prove that half of the structure constants are zero:
\begin{proposition}\label{Prop:vanishing}
For $m$ even, we have $\omega_{k,\ell}^{(m)}=0$.
\end{proposition}

This has been noticed in the proof of Proposition 6.1 in \cite{hitchin1992lie}. 

\begin{proof}
The elements $G_{i,0}$ for $i=1,...,n$ form a basis of $\mathfrak{g}_0=\mathrm{ker}(H)$, which is a Cartan subalgebra of $\mathfrak{sl}_n$. Hence the Lie bracket between $G_{i,0}$ and $G_{j,0}$ has to vanish for all $i,j$. 

Using our procedure from Theorem \ref{Thm:1}, we see that $G_{i,0}$ corresponds to some scalar multiple of $x^iy^i$. Computing the $m$-th transvectant between $x^iy^i$ and $x^jy^j$ gives 0 for $m$ odd and a non-zero multiple of $x^{i+j-m}y^{i+j-m}$ for $m$ even. Since the bracket is 0, the weight of these non-vanishing terms has to be 0.
\end{proof}

Using the procedure, we can give another description of the structure constants by considering a specific Lie bracket between two elements of $V_k$ and $V_\ell$. We choose $E^k\in V_k$ and $F^\ell\in V_\ell$. By definition, they correspond to the monomials
$$(E^k,F^\ell) \cong \left(\frac{x^{2k}}{(2k)!},(-1)^\ell\frac{y^{2\ell}}{(2\ell)!}\right).$$
On the level of polynomials, we know how to compute the transvectant. The projection of $[E^k,F^\ell]$ onto $V_{k+\ell-m}$ is given by
$$\omega_{k,\ell}^{(m)}(-1)^\ell \frac{x^{2k-m}y^{2\ell-m}}{(2k-m)!(2\ell-m)!}\cong\frac{(-1)^\ell \omega_{k,\ell}^{(m)}}{(2\ell-m)!}G_{k+\ell-m,k-\ell}.$$
Therefore, we get
\begin{equation}\label{Eq:useful-eq}
[E^k,F^\ell] = \sum_{m \text{ odd}} \frac{(-1)^\ell \omega_{k,\ell}^{(m)}}{(2\ell-m)!}G_{k+\ell-m,k-\ell}
\end{equation}
which we can also take as definition of the structure constants.

\begin{Remark}
Changing the choice of highest weight vectors from $E^k$ to $\alpha_k E^k$ with $\alpha_k\in \C^*$ for $k=2,...,n-1$, the structure constants change by
$$\omega_{k,\ell}^{(m)}\mapsto \frac{\alpha_k \alpha_\ell}{\alpha_{k+\ell-m}} \omega_{k,\ell}^{(m)}.$$
Hence there are combinations of the structure constants, like $$\frac{\omega_{k,\ell}^{(m)} \omega_{k',\ell'}^{(m')}}{\omega_{k,\ell'}^{(m+\ell'-\ell)} \omega_{k',\ell}^{(m'+\ell-\ell')}}$$ which are independent of the chosen highest weight vectors (for non-vanishing denominator). We have not found interesting structure in these constants though.
\end{Remark}

\section{Trace method}\label{Sec:trace-method}

We present a method which computes some of the structure constants and shows some non-trivial symmetry. The starting point is the following observation:

\begin{proposition}\label{Prop:trace-method}
    Let $(G_{i,j})_{1\leq i\leq n-1, -i\leq j\leq i}$ denote the basis elements of $\mathfrak{sl}_n$ defined in Equation \eqref{def-G}. Then for all $1\leq i,k\leq n-1$ and $-i\leq j\leq i$ and $-k\leq \ell\leq k$ we have $$\tr(G_{i,j}G_{k,\ell}) = 0 \;\text{ if }\; (k,\ell)\neq (i,-j).$$
\end{proposition}

In terms of Figure \ref{Fig:g-decompo}, the proposition says that the trace of a product of two elements of the basis is only non-zero if the two corresponding dots lie symmetric with respect to the middle axis.

\begin{proof}
    The trace gives an isomorphism $V_k\cong V_k^*$ as $\mathfrak{sl}_2$-modules for all $k=1,...,n-1$, via the cyclicity property $\tr ([a,b]c) = \tr([b,c]a)$. Hence if $i\neq k$ and $x\in V_i$, the map $y\mapsto \tr(xy)$ has to be identical zero on $V_k$. This implies the proposition for $i\neq k$.

    Finally, $[H,G_{i,j}G_{i,\ell}]=2(j+\ell)G_{i,j}G_{i,\ell}$, so for $j\neq \ell$ the matrix $G_{i,j}G_{i,\ell}$ is either strictly upper triangular or strictly lower diagonal, in particular of trace zero.
\end{proof}

\begin{proposition}\label{Prop:trace-method-2}
    We have for all $1\leq k\leq n-1$ and $-k\leq \ell\leq k$:
$$\tr(G_{k,\ell}G_{k,-\ell}) = (-1)^{\ell}(2k)! \tr(E^k F^k) = (-1)^{\ell}(2k)! (k!)^2 \binom{n+k}{2k+1}.$$
\end{proposition} 
Note that the result does only depend on the parity of $\ell$.

\begin{proof}
    Using the definition of $G_{k,\ell}$ we get for $-k<\ell<k$:
\begin{align*}
\tr(G_{k,\ell}G_{k,-\ell}) &= \tr(\ad_F(G_{k,\ell+1})G_{k,-\ell})\\
&=-\tr(G_{k,\ell+1}\ad_F(G_{k,-\ell}))\\
&=-\tr(G_{k,\ell+1}G_{k,-\ell-1}).
\end{align*}
Hence $\tr(G_{k,\ell}G_{k,-\ell})=(-1)^{k-\ell}\tr(G_{k,k}G_{k,-k})$. We have $G_{k,k}=E^k$ and by Proposition \ref{prop:lowest-weight-link} also $G_{k,-k}=(-1)^k(2k)!F^k$. The next lemma finishes the proof. 
\end{proof}

\begin{lemma}\label{Lemma-trace}
    We have for $1\leq k\leq n-1$: $$\tr(E^kF^k)=(k!)^2\binom{n+k}{2k+1}.$$
\end{lemma}
The proof is a direct computation using the explicit expression for $F$ and a small combinatorial identity. 
\begin{proof}
    The trace $\tr(E^kF^k)$ is independent of the principal $\mathfrak{sl}_2$-triple since all of them are conjugate. So we can use the one from Equation \eqref{Eq:e-h-f}. Recall that $r_\ell=\ell(n-\ell)$. The only non-zero entries of $F^k$ are $r_jr_{j+1}...r_{j+k-1}$ with $1\leq j\leq n-k$. Hence
    \begin{align*}
        \tr(E^k F^k)&= \sum_{j=1}^{n-k}r_jr_{j+1}...r_{j+k-1} = \sum_{j=1}^{n-k}\frac{(j+k-1)!(n-j)!}{(j-1)!(n-j-k)!}\\
&= (k!)^2\sum_{j=1}^{n-k}\binom{n-j}{k}\binom{j+k-1}{k}.
    \end{align*}
    To compute the remaining sum, imagine you have to choose $2k+1$ objects out of $n+k$ (with $n>k$). Looking at the place $\ell$ of the object number $k+1$, we get
$$\binom{n+k}{2k+1} = \sum_{\ell=k+1}^{n}\binom{\ell-1}{k}\binom{n+k-\ell}{k}=\sum_{j=1}^{n-k}\binom{n-j}{k}\binom{j+k-1}{k}.$$
\end{proof}

Propositions \ref{Prop:trace-method} and \ref{Prop:trace-method-2} allow compute the coefficient of $x\in \mathfrak{sl}_n$ in front of $G_{i,j}$:
$$\mathrm{coeff}_{G_{i,j}}(x)=\frac{\tr(xG_{i,-j})}{\tr(G_{i,j}G_{i,-j})}.$$ 

\begin{proposition}\label{Prop:hidden-symmetry}
    Consider $1\leq \ell\leq k\leq n-1$ and $1\leq m\leq 2\ell-1$ odd. Put $m'=k+\ell-m$. Then, we have the following symmetry for the structure constants:
$$(-1)^{k+1}(2m')!\tr(E^{m'}F^{m'})\omega_{k,\ell}^{(m)} = \tr(E^\ell F^\ell)\frac{(2k-m)!(2\ell)!}{m!}\omega_{m',k}^{(2k-m)}.$$
\end{proposition}
The proof combines Equation \eqref{Eq:useful-eq} with the trace method.
\begin{proof}
    We start from the trace method identity:
\begin{equation}\label{Eq:symmetry}
\tr(G_{m',k-\ell}G_{m',\ell-k})\times \mathrm{coeff}_{G_{m',k-\ell}}([E^k,F^\ell]) = \tr([E^k,F^\ell]G_{m',\ell-k}).
\end{equation}
The coefficient on the left hand side is linked to the structure constant $\omega_{k,l}^{(m)}$ by Equation \eqref{Eq:useful-eq}:
$$\mathrm{coeff}_{G_{m',k-\ell}}([E^k,F^\ell])=\frac{(-1)^\ell}{(2\ell-m)!}\omega_{k,\ell}^{(m)}.$$ In addition, Proposition \ref{Prop:trace-method-2} gives $\tr(G_{m',k-\ell}G_{m',\ell-k})=(-1)^{k+\ell}(2m')!\tr(E^{m'}F^{m'})$. Using the properties of the trace, the right hand side of \eqref{Eq:symmetry} becomes
$$\tr([E^k,F^\ell]G_{m',\ell-k})=\tr(F^\ell[G_{m',\ell-k},E^k])=\tr(E^\ell F^\ell)\mathrm{coeff}_{E^\ell}([G_{m',\ell-k},E^k]).$$
The coefficient of $[G_{m',\ell-k},E^k]$ in front of $E^\ell$ can be computed via our general method. The associated monomials are $\frac{x^{2\ell-m}y^{2k-m}}{(2\ell-m)!}$ and $\frac{x^{2k}}{(2k)!}$. Applying the transvectant of order $2k-m$ gives
$$(-1)^{2k-m}\omega_{m',k}^{(2k-m)}\frac{(2k-m)!}{(2\ell-m)!m!}x^{2\ell}.$$
Hence $\mathrm{coeff}_{E^\ell}([G_{m',\ell-k},E^k])=-\omega_{m',k}^{(2k-m)}\frac{(2k-m)!(2\ell)!}{(2\ell-m)!m!}$. The sign comes from the fact that $m$ is odd.
Combining all together, we get
$$(-1)^{k+1}\frac{(2m')!}{(2\ell-m)!}\tr(E^{m'}F^{m'})\omega_{k,\ell}^{(m)} = \tr(E^\ell F^\ell)\frac{(2k-m)!(2\ell)!}{(2\ell-m)!m!}\omega_{m',k}^{(2k-m)}$$
which proves the proposition.
\end{proof}

Using the trace method, it is possible to compute the structure constants $\omega_{k,\ell}^{(m)}$ for the extremal values of $m$, i.e. $m=1$ and $m=2\ell-1$. It seems impossible to get a general formula via this method. This is why we use the more powerful graphical calculus.

\section{Graphical calculus}\label{Sec:graphical-calculus}

\subsection{Introduction to the graphical calculus for classical invariant theory}

The following explicit computations will use the graphical formalism developed in~\cite[\S2 and \S3]{AbdesselamJKTR} (see also~\cite[\S4.2]{AbdesselamCDCG} for additional general explanations) with the purpose of putting under the same roof the 19th century computational techniques from the classical invariant theory of binary forms, and in particular the classical symbolic method, as well as the 20th century theory of quantum angular momentum in mathematical physics.
We will work with very concrete ``old-fashioned'' tensors simply seen as multidimensional arrays of numbers with indices belonging to the two element set $[2]:=\{1,2\}$.
We will build more complicated tensors out of some basic building blocks, using contraction of indices, but since such expressions quickly become unwieldy, we will use diagrams as a numerically precise shorthand notation for these expressions.

We will denote points in $\bb{C}^2$ or pairs of variables by lowercase
boldface letters such as $\mathbf{x}=(x_1,x_2)$, $\mathbf{y}=(y_1,y_2)$.
A \emph{binary form} of degree $d$ is a polynomial function $F(\mathbf{x})$ of $x_1,x_2$ which is homogeneous of degree $d$.
We will identify such a function with the uniquely defined fully symmetric tensor
\[
F=(F_{i_1,\ldots,i_d})_{i_1,\ldots,i_d\in[2]}
\]
which satisfies
\[
F(\mathbf{x})=\sum_{i_1,\ldots,i_d\in[2]}
F_{i_1,\ldots,i_d}\ x_{i_1}\cdots x_{i_d}\ ,
\]
for all $\mathbf{x}\in\bb{C}^2$.
By fully symmetric tensor or array we mean that $F$ satisfies
\[
F_{i_{\sigma(1)},\ldots,i_{\sigma(d)}}=F_{i_1,\ldots,i_d}\ ,
\]
for all permutations $\sigma\in\mathfrak{S}_d$ and for all specifications $(i_1,\ldots,i_d)\in[2]^d$ of the indices.
As in~\cite{AbdesselamJKTR}, we denote by $\mc{H}_{d}$ the space of binary forms of degree $d$ and define the $\SL_2(\bb{C})$ action on it as follows.
For a matrix
\[
g=\begin{pmatrix}g_{11} & g_{12}\\ g_{21} & g_{22}\end{pmatrix}
\]
in $\SL_2(\C)$, we define its action on vectors $\mathbf{x}=(x_1,x_2)$
by 
\[
g\mathbf{x}:=(g_{11}x_1+g_{12}x_2,g_{21}x_1+g_{22}x_2)\ ,
\]  
namely, $(g\mathbf{x}^{\rm T})^{\rm T}$, in terms of matrix algebra. In other words, we think of $\mathbf{x}$ as a column vector when computing matrix products, but we will write it as a row vector, especially when it is the argument of a polynomial function $F$.
We then define the action on binary forms by letting
\[
(gF)(\mathbf{x}):=F(g^{-1}\mathbf{x})
\]
for all $\mathbf{x}$.
As is well known, $\mc{H}_{d}$ with this action is a concrete model 
for the $(d+1)$-dimensional irreducible representation of $\SL_2(\C)$. Comparing to the notation in Section \ref{Sec:natural-basis}, we get $V_d\cong \mc{H}_{2d}$.
Using the identification with tensors, this action can also be defined by
\begin{equation}
(gF)_{i_1,\ldots,i_d}=\sum_{j_1,\ldots,j_d\in[2]}
(g^{-1})_{j_1 i_1}\cdots(g^{-1})_{j_d i_d}
\ F_{j_1,\ldots,j_d}
\label{Ftransclassicaleq}
\end{equation}
which, in the graphical formalism of~\cite{AbdesselamJKTR}, becomes the equation
\begin{equation}
\parbox{2.6cm}{
\includegraphics[width=2.6cm]{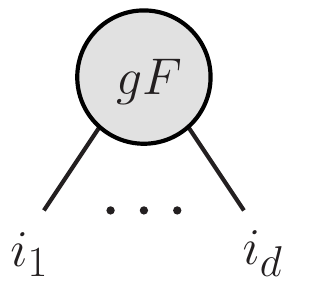}}
 := 
\parbox{4.6cm}{
\includegraphics[width=4.6cm]{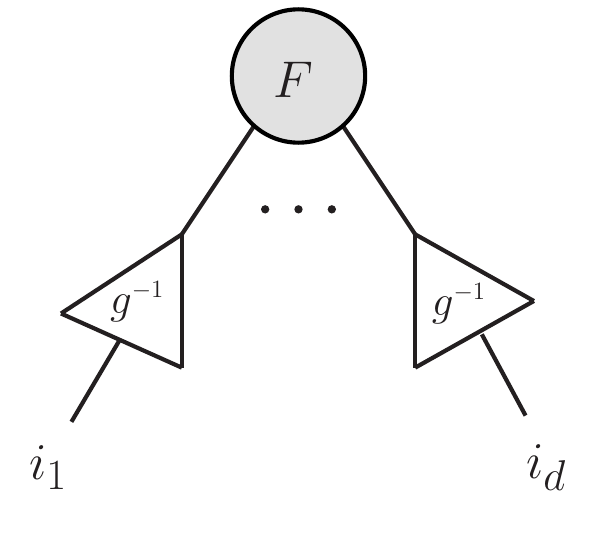}}
\ .
\label{Ftransgraphicaleq}
\end{equation}
For a binary form $F$ of degree $d$, we use a round ``blob'' with $d$ ``legs'' attached in order to denote the tensor entry $F_{i_1,\ldots,i_d}$. Namely,
\[
\parbox{2.6cm}{
\includegraphics[width=2.6cm]{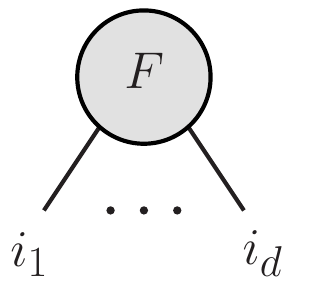}}
:= F_{i_1,\ldots,i_d}\ .
\]
For a $2\times 2$ matrix $g=(g_{ij})_{i,j\in[2]}$, we introduce the following graphical notation
\[
\parbox{2.8cm}{\raisebox{-11ex}{
\includegraphics[width=2.8cm]{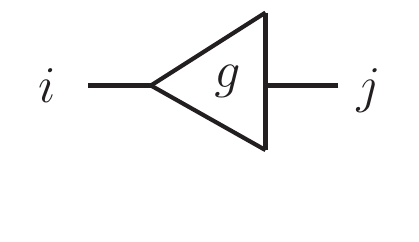}}
}
:=g_{ij}
\] 
for its entries.
For the identity matrix $g=I$, we just put a line with no triangular box
\[
\parbox{2.5cm}{\raisebox{-10ex}{
\includegraphics[width=2.5cm]{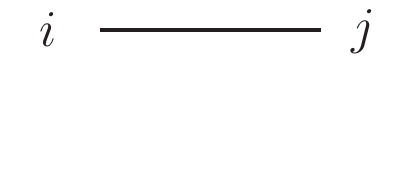}}
}
:=\delta_{ij}\ ,
\]
using the standard Kronecker delta symbol.
For the important special matrix
\[
\varepsilon=\begin{pmatrix}
0 & 1\\
-1 & 0
\end{pmatrix}\ ,
\]
we use an arrow instead of a triangular box
\[
\parbox{2.5cm}{\raisebox{-10ex}{
\includegraphics[width=2.5cm]{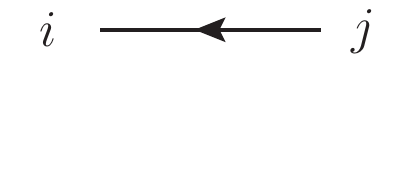}}
}
:=\varepsilon_{ij}\ .
\]
The above lists the basic building blocks of the graphical calculus. The way to evaluate more complicated pictures made of such blocks is to first take the product of the corresponding matrix or tensor entries for the blocks present, and then, whenever two
legs are glued, one is to assign the same index to both glued legs and sum over the two possible values of that index in the set $[2]$. One must do that independently, for every pair of glued legs.
For example, if $\mathbf{y}=g\mathbf{x}$, then
saying
\[
\parbox{1.7cm}{\raisebox{-11ex}{
\includegraphics[width=1.7cm]{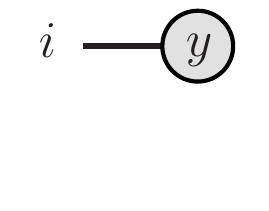}}
}
=\parbox{2.8cm}{\raisebox{-11ex}{
\includegraphics[width=2.8cm]{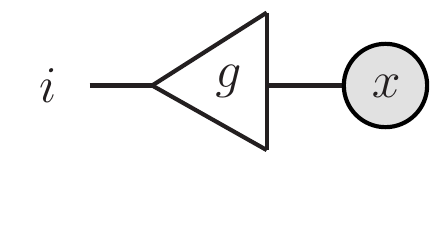}}
}
\ ,
\]
for all values of the free leg index $i\in[2]$, is the same as writing
\[
y_i=\sum_{j\in[2]}g_{ij}x_j\ ,
\]
which defines the transformed vector $\mathbf{y}$.
In the following computations, a fundamental role is played by another basic building block which is the \emph{symmetrizer} of size $n$
\[
\parbox{1.7cm}{\raisebox{-11ex}{
\includegraphics[width=1.7cm]{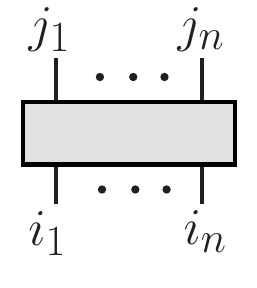}}
}
\ \ :=\mc{S}_{i_1,\ldots,i_n}^{j_1,\ldots,j_n}
:=\frac{1}{n!}\sum_{\sigma\in\mf{S}_n}\delta_{i_1,j_{\sigma(1)}}\cdots\delta_{i_n,j_{\sigma(n)}}
\]
where, as usual, $\mf{S}_n$ is the symmetric group on $n$ elements.

A key notion of classical invariant theory is that of \emph{transvectant} of order $k$ of two binary forms $F,G$ of respective degrees $m,n$. The allowed range for $k$ is $0\le k\le \min(m,n)$.
As a polynomial in $\mathbf{x}$, and with the classical normalization, it is defined by
\begin{equation}\label{Eq:transvectant-graphical}
(F,G)_k:=
\parbox{4.3cm}{\raisebox{-11ex}{
\includegraphics[width=4.3cm]{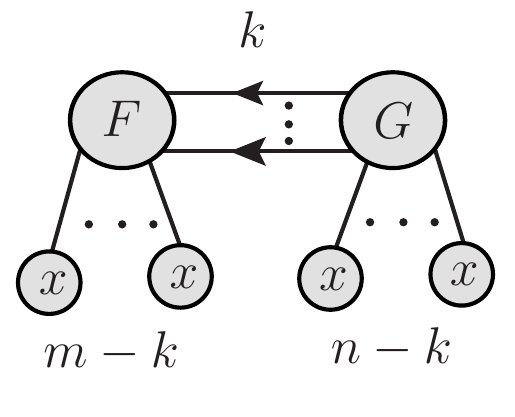}}
}
\ .
\end{equation}
The above graphical shorthand equation, by itself, is a mathematically precise definition. However, for the reader who is not familiar with this kind of diagrammatic formulas, let us give, as a ``cheat sheet'', the longhand form of the same statement
\[
(F,G)_k:=\sum_{i_1,\ldots,i_m,j_1,\ldots,j_n\in[2]}F_{i_1,\ldots,i_m}\ \varepsilon_{i_1,j_1}\cdots\varepsilon_{i_k,j_k}\ G_{j_1,\ldots,j_n}\ 
x_{i_{k+1}}\cdots x_{i_m}\ x_{j_{k+1}}\cdots x_{j_n}\ .
\]
Another formula for this transvectant, as a differential operator acting on a polynomial expression, is
\[
(F,G)_k=
\frac{(m-k)!\ (n-k)!}{m!\ n!}
\left.(\Omega_{\mathbf{x}\mathbf{y}})^k F(\mathbf{x})G(\mathbf{y})
\right|_{\mathbf{y}:=\mathbf{x}}
\]
which features Cayley's Omega Operator
\[
\Omega_{\mathbf{x}\mathbf{y}}:=\frac{\partial^2}{\partial x_1 \partial y_2}-
\frac{\partial^2}{\partial x_2 \partial y_1}\ .
\]
Still another useful formula, in particular for computer implementation, is
\begin{equation}
(F,G)_k(\mathbf{x})=
\frac{(m-k)!\ (n-k)!}{m!\ n!}\ 
\sum_{j=0}^k
(-1)^j
\binom{k}{j}
\frac{\partial^k F}{\partial x_1^{k-j}\partial x_2^{j}}(\mathbf{x})
\frac{\partial^k G}{\partial x_1^{j}\partial x_2^{k-j}}(\mathbf{x})\ .
\label{transveq}
\end{equation}
It is not hard to see that transvectants satisfy, for all $g\in\SL_2$, and all binary forms $F$, $G$ of the specified degrees, the following identity
\[
g\cdot (F,G)_k=(gF, gG)_k\ ,
\]
namely, the construction is $\SL_2$-equivariant, and explicitly realizes the projection from $\mc{H}_{m}\otimes\mc{H}_{n}$ onto the irreducible representation $\mc{H}_{m+n-2k}$.
In fact, we will make this projection even more explicit by thinking of the modern tensor product $\mc{H}_{m}\otimes\mc{H}_{n}$ as the more concrete space of bihomogeneous polynomials $A(\mathbf{x},\mathbf{y})$
in two sets of variables $\mathbf{x}$, $\mathbf{y}$ of bidegree $(m,n)$.
We will identify such a bihomogeneous form $A$ with its (old-fashioned) tensor
\[
A=(A_{i_1,\ldots,i_m;j_1,\ldots,j_n})_{i_1,\ldots,i_m,j_1,\ldots,j_n\in[2]}
\]
which must be $\mf{S}_{m}\times\mf{S}_{n}$-symmetric, i.e., must be satisfy
\[
A_{i_{\sigma(1)},\ldots,i_{\sigma(m)};j_{\tau(1)},\ldots,j_{\tau(n)}}
=A_{i_1,\ldots,i_m;j_1,\ldots,j_n}
\]
for all permutations $\sigma\in\mf{S}_m$ and $\tau\in\mf{S}_n$, and all values of the $i$ and $j$ indices in the set $[2]$. The identification between the $\mf{S}_{m}\times\mf{S}_{n}$-symmetric array $A$ and the corresponding bihomogeneous form $A(\mathbf{x},\mathbf{y})$ is given by imposing that for all $\mathbf{x}=(x_1,x_2)$ and $\mathbf{y}=(y_1,y_2)$ in $\mathbb{C}^2$, we have
\[
A(\mathbf{x},\mathbf{y})=\sum_{(i_1,\ldots,i_m,j_1,\ldots,j_n)\in[2]^{m+n}}
A_{i_1,\ldots,i_m;j_1,\ldots,j_n}\ 
x_{i_1}\cdots x_{i_m}\ y_{j_1}\cdots y_{j_n}\ .
\]
We will also introduce a ``SIM card'' notation for the corresponding tensor entries, namely,
\begin{equation}
\parbox{2.5cm}{\raisebox{-11ex}{
\includegraphics[width=2.5cm]{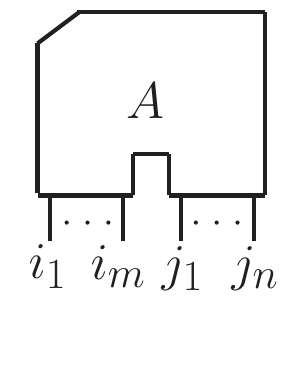}}
}
\ \ 
:= A_{i_1,\ldots,i_m;j_1,\ldots,j_n}\ .
\label{simcardeq}
\end{equation}

In the upcoming calculations,
we will use two fundamental graphical identities.
The first one
\begin{equation}
\parbox{4.7cm}{\raisebox{-11ex}{
\includegraphics[width=4.7cm]{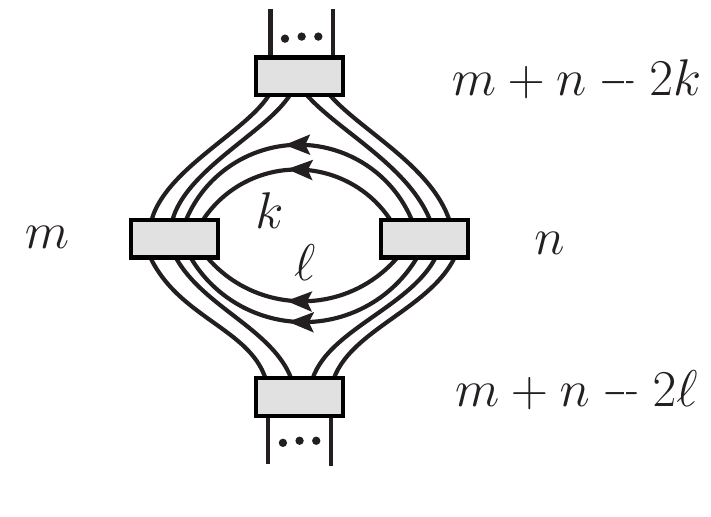}}
}
=\delta_{k,\ell}\times
\frac{\binom{m+n-k+1}{k}}{\binom{m}{k} \binom{n}{k}}\times
\parbox{3.7cm}{\raisebox{-11ex}{
\includegraphics[width=3.7cm]{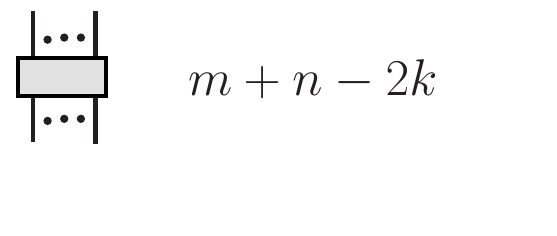}}
}
\label{CGSchureq}
\end{equation}
is a diagrammatic version of Schur's Lemma, see~\cite[Eq. (2.10)]{AbdesselamJKTR}.
The second one
\begin{equation}
\parbox{3.3cm}{\raisebox{-12ex}{
\includegraphics[width=3.3cm]{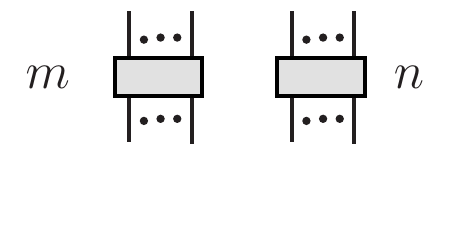}}
}
=\sum_{k=0}^{\min(m,n)}
\frac{\binom{m}{k} \binom{n}{k}}{\binom{m+n-k+1}{k}}
\times
\parbox{4.7cm}{\raisebox{-11ex}{
\includegraphics[width=4.7cm]{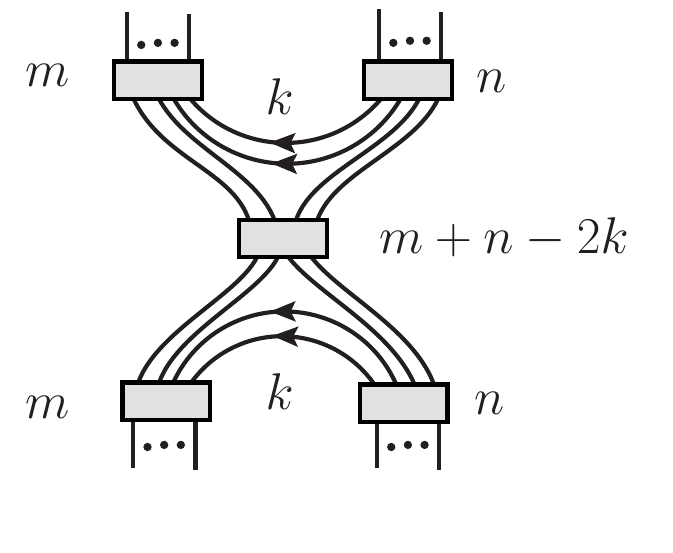}}
}
\label{CGdecompeq}
\end{equation}
is an explicit form of the Clebsch-Gordan decomposition of tensor products of irreducible $\SL_2$-representations, see~\cite[Eq. (2.9)]{AbdesselamJKTR}. Namely, the decomposition is carried out using explicit intertwiners.
Let us again provide the nongraphical ``cheat sheet'' for the graphical equation (\ref{CGdecompeq}). The precise statement which the latter says is that $\forall(a_1,\ldots,a_m)\in[2]^m$, $\forall(b_1,\ldots,b_n)\in[2]^n$, $\forall(c_1,\ldots,c_m)\in[2]^m$, $\forall(d_1,\ldots,d_m)\in[2]^n$,
\[
L_{c_1,\ldots,c_m;d_1,\ldots,d_n}^{a_1,\ldots,a_m;b_1,\ldots,b_n}=R_{c_1,\ldots,c_m;d_1,\ldots,d_n}^{a_1,\ldots,a_m;b_1,\ldots,b_n}\ ,
\]
where the left-hand side is, by definition,
\[
L_{c_1,\ldots,c_m;d_1,\ldots,d_n}^{a_1,\ldots,a_m;b_1,\ldots,b_n}:=
\mathcal{S}_{c_1,\ldots,c_m}^{a_1,\ldots,a_m}\ \times\ \mathcal{S}_{d_1,\ldots,d_n}^{b_1,\ldots,b_n}\ ,
\]
while the right-hand side is, by definition,
\begin{eqnarray*}
    R_{c_1,\ldots,c_m;d_1,\ldots,d_n}^{a_1,\ldots,a_m;b_1,\ldots,b_n} & := &
    \sum_{k=0}^{\min(m,n)}
\frac{\binom{m}{k} \binom{n}{k}}{\binom{m+n-k+1}{k}}
\times \\
& & \sum_{(e_1,\ldots,e_m)\in[2]^m}
\ \sum_{(f_1,\ldots,f_n)\in[2]^n}
\ \sum_{(g_1,\ldots,g_m)\in[2]^m}
\ \sum_{(h_1,\ldots,h_n)\in[2]^m} \\
 & & \ \ \ \mathcal{S}_{e_1,\ldots,e_m}^{a_1,\ldots,a_m}
\ \mathcal{S}_{f_1,\ldots,f_n}^{b_1,\ldots,b_n} \\
 & & \times\ \varepsilon_{e_{m-k+1},f_{k}}\ \varepsilon_{e_{m-k+2},f_{k-1}}\cdots
\varepsilon_{e_m,f_1} \\
& & \times\ \mathcal{S}_{g_1,\ldots,g_{m-k},h_{k+1},\ldots,h_n}^{e_1,\ldots,e_{m-k},f_{k+1},\ldots,f_n} \\
 & & \times\ \varepsilon_{g_{m-k+1},h_{k}}\ \varepsilon_{g_{m-k+2},h_{k-1}}\cdots
\varepsilon_{g_m,h_1} \\
 & & \times\ \mathcal{S}_{c_1,\ldots,c_m}^{g_1,\ldots,g_m}
\ \mathcal{S}_{d_1,\ldots,d_n}^{h_1,\ldots,h_n}\ .
\end{eqnarray*}
Thus the graphical identity (\ref{CGdecompeq}) packages together a collection of $2^{2m+2n}$ equalities between scalars, and this collection can be interpreted as asserting the equalities of two maps. Indeed, one can use the left-hand side to define a map which sends the $\mf{S}_m\times\mf{S}_n$-symmetric tensor or array $A$ to a new one $B$ defined by 
letting $\forall(c_1,\ldots,c_m)\in[2]^m$, and $\forall(d_1,\ldots,d_n)\in[2]^n$,
\[
B_{c_1,\ldots,c_m;d_1,\ldots,d_n}:=
\sum_{(a_1,\ldots,a_m,b_1,\ldots,b_n)\in[2]^{m+n}}
L_{c_1,\ldots,c_m;d_1,\ldots,d_n}^{a_1,\ldots,a_m;b_1,\ldots,b_n}\ 
A_{a_1,\ldots,a_m;b_1,\ldots,b_n}\ .
\]
One can define a similar map using $R_{c_1,\ldots,c_m;d_1,\ldots,d_n}^{a_1,\ldots,a_m;b_1,\ldots,b_n}$ instead of $L_{c_1,\ldots,c_m;d_1,\ldots,d_n}^{a_1,\ldots,a_m;b_1,\ldots,b_n}$. In this instance, the map on the left-hand side is the identity on $\mathcal{H}_{m}\otimes\mathcal{H}_{n}$, whereas the map on the right-hand side is given as a linear combination of maps of the form $\iota\circ\pi$ where $\pi$ is an ${\rm SL}_2$-equivariant projection $\mathcal{H}_{m}\otimes\mathcal{H}_{n}\rightarrow\mathcal{H}_{m+n-2k}$ (the top half of the picture) and $\iota$ is an ${\rm SL}_2$-equivariant injection $\mathcal{H}_{m+n-2k}\rightarrow\mathcal{H}_{m}\otimes\mathcal{H}_{n}$ (the bottom half of the picture).
For a proof of equivariance of these maps see \cite[Prop. 3.1]{AbdesselamJKTR}.

Finally, we will also work with ${\rm End}(\mc{H}_d)=\mc{H}_d\otimes\mc{H}_d^{\vee}$, the space of linear maps from the space of binary forms of degree $d$ to itself. Such a map $M$ will also be identified
with a tensor $M=(M_{i_1,\ldots,i_d;j_1,\ldots,j_d})_{i_1,\ldots,i_d,j_1,\ldots,j_d\in[2]}$ with the imposed $\mf{S}_d\times\mf{S}_d$ symmetry of invariance by permutation of the $i$ indices, and by permutation of the $j$ indices.
The relation between the map and the tensor, in graphical terms is that for all binary forms $F$ of degree $d$, we have
\[
\parbox{2.5cm}{\raisebox{-20ex}{
\includegraphics[width=2.5cm]{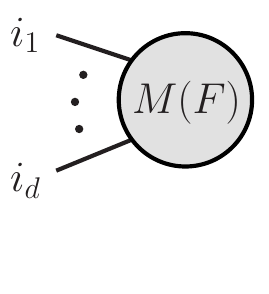}}
}
\ \ := \ 
\parbox{4.5cm}{\raisebox{-17ex}{
\includegraphics[width=4.5cm]{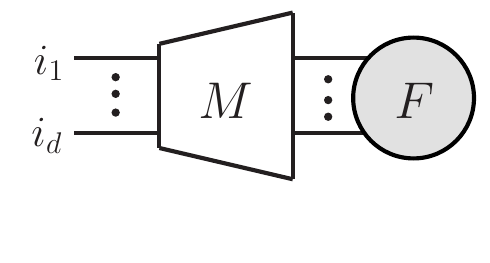}}
}
\ \ := \sum_{j_1,\ldots,j_d\in[2]}
M_{i_1,\ldots,i_d;j_1,\ldots,j_d} \ F_{j_1,\ldots,j_d}\ ,
\]
for all values of the indices $i_1,\ldots,i_d$ in $[2]$.
Here we introduced another graphical piece of notation, with a trapeze-like shape, for the $M$ tensor, namely,
\[
\parbox{4.2cm}{\raisebox{-17ex}{
\includegraphics[width=4.2cm]{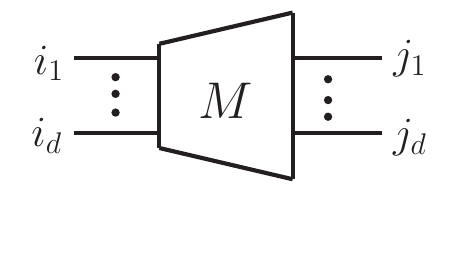}}
}
:= 
M_{i_1,\ldots,i_d;j_1,\ldots,j_d}\ .
\]
It is not difficult to see that the trace of the endomorphism $M$ is then given by
\[
\tr(M)= 
\parbox{4cm}{\raisebox{-11ex}{
\includegraphics[width=4cm]{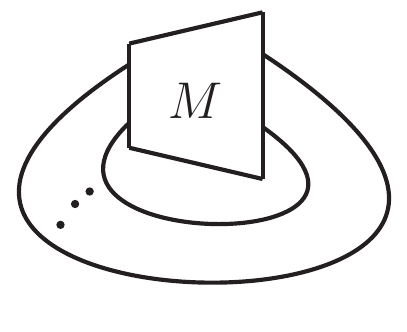}}
}
\ :=\sum_{i_1,\ldots,i_d\in[2]}
M_{i_1,\ldots,i_d;i_1,\ldots,i_d} \ .
\]
Note that composition of endomorphisms $M\circ N$ can easily be expressed graphically  by
\[
\parbox{3cm}{\raisebox{-16ex}{
\includegraphics[width=3cm]{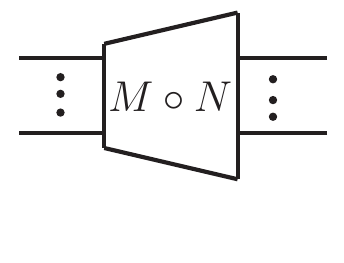}}
}
\ \ =\ \  
\parbox{4.7cm}{\raisebox{-15ex}{
\includegraphics[width=4.7cm]{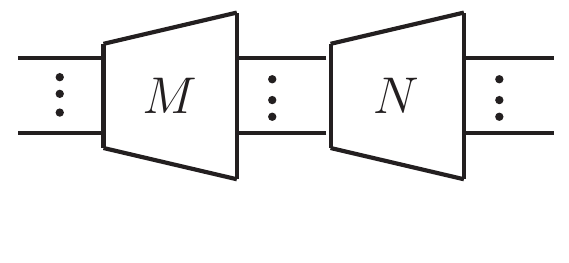}}
}
\ \ .
\]

We will take advantage of self-duality of $\SL_2$-representations in an explicit manner, as follows. We will identify an endomorphism $M$ to the bihomogeneous form $A$ of bidegree $(d,d)$ defined by the diagrammatic equation
\begin{equation}
\parbox{2.5cm}{\raisebox{-11ex}{
\includegraphics[width=2.5cm]{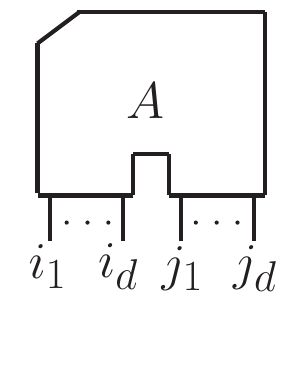}}
}
\ \ := 
\parbox{5cm}{\raisebox{-11ex}{
\includegraphics[width=5cm]{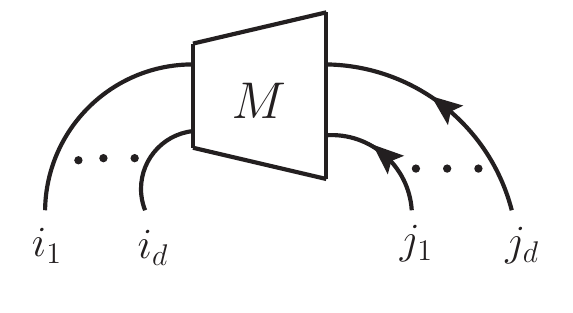}}
}
\ ,
\label{MtoAeq}
\end{equation}
i.e.,
\[
A_{i_1,\ldots,i_d;j_1,\ldots,j_d} =\sum_{\ell_1,\ldots,\ell_d\in[2]}
M_{i_1,\ldots,i_d;\ell_1,\ldots,\ell_d}\ 
\varepsilon_{\ell_1 j_1}\cdots\varepsilon_{\ell_d j_d}\ , 
\]
for all values of the $i$ and $j$ indices.
As a result, and abusing notation by writing $A=M$ etc., the trace operation, in the bihomogeneous form point-of-view becomes
\begin{equation}
\tr(A)= 
\parbox{2.4cm}{\raisebox{-11ex}{
\includegraphics[width=2.4cm]{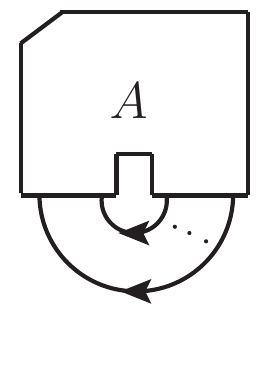}}
}
\ \ \ ,
\label{simtraceeq}
\end{equation}
and composition becomes
\begin{equation}
\parbox{2.5cm}{\raisebox{-13ex}{
\includegraphics[width=2.5cm]{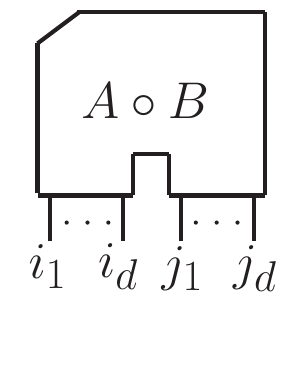}}
}
\ \ =\ \  
\parbox{6.2cm}{\raisebox{-11ex}{
\includegraphics[width=6.3cm]{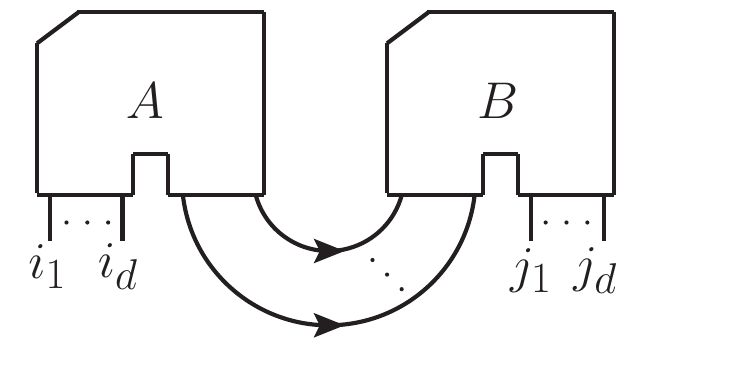}}
}
\ \ \ .
\label{simcompeq}
\end{equation}

\subsection{Graphical calculus for Wigner's $6j$ symbols}

The \emph{Wigner $6j$ symbol} denoted by
\[
\left\{
\begin{array}{ccc}
j_1 & j_2 & j_{12}\\
j_3 & J & j_{23}
\end{array}
\right\}
\]
encodes a precise standard numerical function of the 6 entries $j_1,j_2,j_{12},j_3,J,j_{23}$ which belong to $\frac{1}{2}\bb{N}$.
The domain is defined by the requirement that the four triples $(j_1,j_2,j_{12})$, $(j_2,j_3,j_{23})$, $(j_1,j_{23},J)$, $(j_{12},j_3,J)$
must be triads. We say that a triple such as $(j_1,j_2,j_{12})$
is a \emph{triad} iff $j_1+j_2+j_{12}\in\bb{N}$
and $|j_1-j_2|\le j_{12}\le j_1+j_2$.
The standard definition of $6j$ symbols, together with the explanation and motivation for the choices of conventions (e.g., the Condon-Shortley phase convention)
involved in this standard definition, are recalled in~\cite[\S7]{AbdesselamCAIF}.
The relation to our graphical calculus is explained in~\cite[\S7.1 and \S7.2]{AbdesselamCMMJlong}.

Consider the map $\psi:\mc{H}_{2J}\rightarrow\mc{H}_{2J}$ which sends a binary form $F$ of degree $2J$ to the binary form $G$ with associated symmetric tensor graphically defined by
\begin{equation}
G_{i_1,\ldots,i_{2J}}:= \ \ 
\parbox{7.5cm}{\raisebox{-11ex}{
\includegraphics[width=7.5cm]{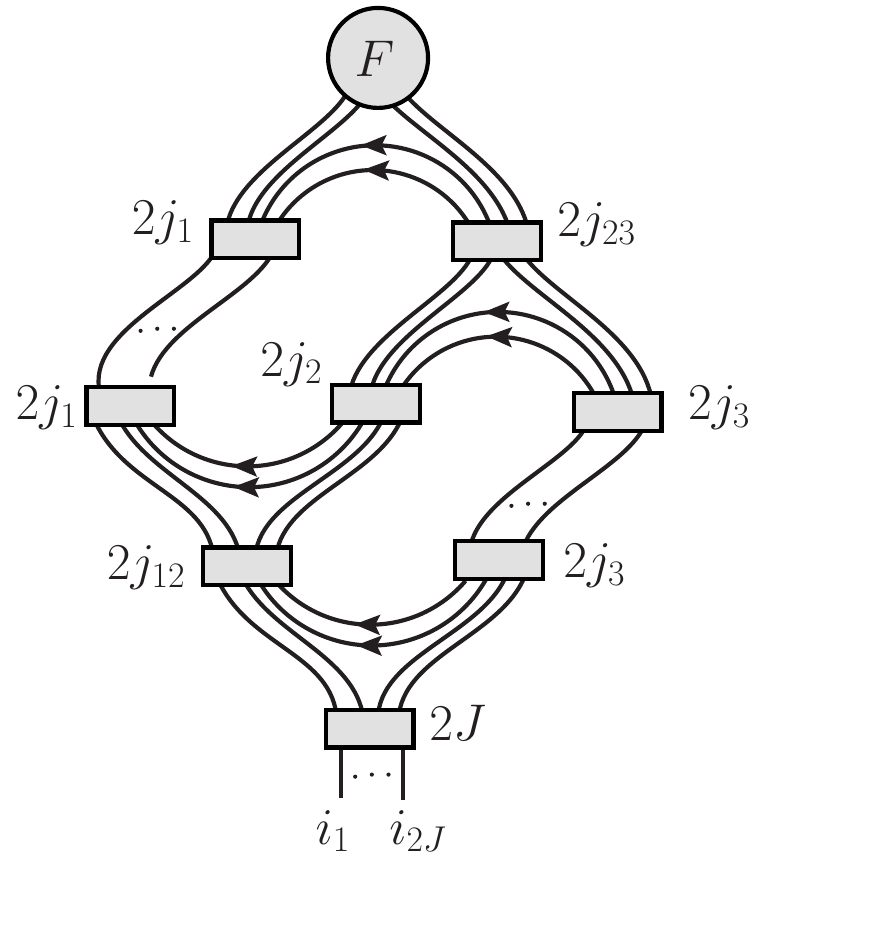}}
}
\ \ \ ,
\label{6jdefeq}
\end{equation}
for all values of the indices $i_1,\ldots,i_{2J}$ in $[2]=\{1,2\}$.
The numbers $2j_1$, etc. on the picture indicate the size of the symmetrizer (number of strands passing through it) next to that number. We did not write the number of $\varepsilon$ arrows so as not to overload the picture. These numbers of arrows are uniquely determined by a simple counting of strands, and are given, from top to bottom by $j_1+j_{23}-J$, $j_2+j_3-j_{23}$, $j_1+j_2-j_{12}$, and $j_{12}+j_3-J$.
The map $\psi$ is $\SL_2$-equivariant and therefore is equal to a multiple of the identity, i.e., $\psi=\rho\times Id$ for a specific value of the constant 
\[
\rho=:\rho\left[
\begin{array}{ccc}
j_1 & j_2 & j_{12}\\
j_3 & J & j_{23}
\end{array}
\right]\ .
\]
By definition, the standard $6j$ symbol is equal to
\[
\left\{
\begin{array}{ccc}
j_1 & j_2 & j_{12}\\
j_3 & J & j_{23}
\end{array}
\right\}
:=(-1)^{j_1+j_2+j_3+J}
(2J+1)\frac{K_1}{\sqrt{K_2 K_3}}\times
\rho\left[
\begin{array}{ccc}
j_1 & j_2 & j_{12}\\
j_3 & J & j_{23}
\end{array}
\right]\ ,
\]
where
\[
K_1:=(2j_1)!(2j_2)!(2j_3)!(2j_{12})!(2j_{23})!(2J)!\ ,
\]
\[
K_2:=(j_1+j_2+j_{12}+1)!(j_2+j_3+j_{23}+1)!(j_1+j_{23}+J+1)!
(j_{12}+j_3+J+1)!\ ,
\]
\begin{eqnarray*}
K_3:= &  & (j_1+j_2-j_{12})!(j_{1}+j_{12}-j_{2})!(j_{2}+j_{12}-j_{1})!\\
 & \times & (j_2+j_3-j_{23})!(j_{2}+j_{23}-j_{3})!(j_{3}+j_{23}-j_{2})!\\
 & \times & (j_1+j_{23}-J)!(j_{1}+J-j_{23})!(j_{23}+J-j_{1})!\\
 & \times & (j_{12}+j_3-J)!(j_{12}+J-j_{3})!(j_{3}+J-j_{12})!\ .
\end{eqnarray*}

The standard $6j$ symbol enjoys a large number of symmetries, and this is a reason why one includes odd-looking factors such as $(-1)^{j_1+j_2+j_3+J}
(2J+1)$ in the definition. The $6j$ symbol is invariant by any permutation of the columns. It is also invariant by simultaneously flipping the top and bottom entries in any two columns. There are also other more complicated Regge symmetries,
but we will not need them.

The simplest known formula for these symbols is Racah's celebrated single sum formula as a terminating ${}_4F_3$ hypergeometric series:
\begin{equation}
\left\{
\begin{array}{ccc}
\! j_1 & j_2 & j_{12} \! \\
\! j_3 & J & j_{23} \!
\end{array}
\right\}
=\sqrt{\frac{K_3}{K_2}}\sum_q\frac{(-1)^q\ (q+1)!}{(q\! -\! T_1)!(q\! -\! T_2)!(q\! -\! T_3)!(q\! -\! T_4)!
(S_1\! -\! q)!(S_2\! -\! q)!(S_3\! -\! q)!}\ ,
\label{6jformulaeq}
\end{equation}
where
\begin{eqnarray*}
T_1 & := & j_1+j_2+j_{12}\ ,\\
T_2 & := & j_2+j_3+j_{23}\ ,\\
T_3 & := & j_1+j_{23}+J\ ,\\
T_4 & := & j_{12}+j_3+J\ ,
\end{eqnarray*}
and
\begin{eqnarray*}
S_1 & := & j_1+j_2+j_3+J\ ,\\
S_2 & := & j_2+j_{12}+j_{23}+J\ ,\\
S_3 & := & j_1+j_3+j_{12}+j_{23}\ .
\end{eqnarray*}

The range of summation is $q\in\bb{Z}$ together with the requirement that the arguments of all the factorials are nonnegative. Note that the condition $q+1\ge 0$ however is redundant, because the $T$'s are nonnegative.

For our Lie algebra computations, we will mostly need the quantity $\rho$ which therefore is given by the formula
\[
\rho\left[
\begin{array}{ccc}
j_1 & j_2 & j_{12}\\
j_3 & J & j_{23}
\end{array}
\right]
=
\frac{(-1)^{j_1+j_2+j_3+J}K_3}{(2J+1)\ K_1}
\]
\begin{equation}
\times
\sum_q\frac{(-1)^q (q+1)!}{(q-T_1)!(q-T_2)!(q-T_3)!(q-T_4)!
(S_1-q)!(S_2-q)!(S_3-q)!}\ .
\label{rhoeq}
\end{equation}

\smallskip
The proof of the above single-sum formula for the $6j$ symbol, due to Racah, is nontrivial and will not be given here. For the diligent reader who would like too see a proof of this formula, we suggest two approaches.

\smallskip
\noindent{\bf 1st approach: Racah's original proof.}
The graphical definition shows that the $6j$ symbol is a contraction of four $3jm$ symbols (see~\cite[\S7.5]{AbdesselamCAIF} for their standard definition). These are the trilinear objects corresponding to transvectants seen in the basis of monomials in $x_1,x_2$. They are given by a single-sum formula, see~\cite[\S5]{AbdesselamCTG}.
Using these ingredients, one can then follow the proof given by Racah in~\cite[App. B]{Racah} and which uses the Chu-Vandermonde Theorem and variants, in both directions.

\smallskip
\noindent{\bf 2nd approach: The Penrose-Moussouris chromatic method.}
This approach was briefly described in~\cite[\S7.2]{AbdesselamCMMJlong} (not to be confused with the shorter published version~\cite{AbdesselamCMMJshort}).
One starts by using~\cite[Thm. 4.1]{AbdesselamJKTR}
which expresses the $6j$ symbol in terms of Penrose's spin network evaluation. In the latter, symmetrizers become antisymmetrizers and loops, i.e., traces of the identity contribute a factor $-2$ instead of $2$.
In~\cite[Thm. 4.1]{AbdesselamJKTR}, there is a global sign which has a complicated expression, for general spin networks, but is easy to compute for the particular case of the $6j$ symbol (see~\cite[\S7.2]{AbdesselamCMMJlong}). Then the formula reduces to Lemma 2.3 of~\cite{GaroufalidisvdVZ}. The proof of this lemma, following the method of Penrose and Moussouris, is given in~\cite[Prop. 12]{KauffmanL}, with the specialization $N=-2$ for the number of ``colors''.

\subsection{The graphical calculus implementation of the principal $\mf{sl}_2$ triple}\label{Sec:graphical-calc-for-sln}

We pick up the thread from Section \ref{Sec:structure-constants}, and give an explicit realization of the Lie algebras $\mf{sl}_n$, and $\mf{gl}_n$, with $n\ge 2$, in terms of the invariant theory of binary forms and the previous graphical calculus.
We will see $\mf{sl}_n$ as sitting inside $\mf{gl}_n$ while viewing the latter as ${\rm End}(\mc{H}_{n-1})$, as a first step. However, as a second step, we will immediately switch to the bihomogeneous form point-of-view and identify $\mf{gl}_n$ with $\mc{H}_{n-1}\otimes\mc{H}_{n-1}$ instead of $\mc{H}_{n-1}\otimes\mc{H}_{n-1}^{\vee}$. Thus, if we say $A\in\mf{gl}_n$, we mean that $A=A(\mathbf{x},\mathbf{y})$ is a bihomogeneous form of bidegree $(n-1,n-1)$.
The element $A$ will therefore have an associated tensor and graphical ``SIM card'' representation as in (\ref{simcardeq}).
The Lie bracket of $\mf{gl}_n$
is of course given by
\[
[A,B]=A\circ B-B\circ A\ ,
\]
where the composition of endomorphisms, now seen as bihomogeneous forms, is computed using formula (\ref{simcompeq}).
We now introduce a map $\Pi:\mf{gl}_n\rightarrow \mc{H}_{0}\times\mc{H}_{2}\times\mc{H}_{4}\cdots\times\mc{H}_{2n-2}$ which send $A$ to $\Pi(A)=(A_0,A_1,\ldots,A_{n-1})$ where
\[
A_i:= \ \ 
\parbox{3cm}{\raisebox{-12ex}{
\includegraphics[width=3cm]{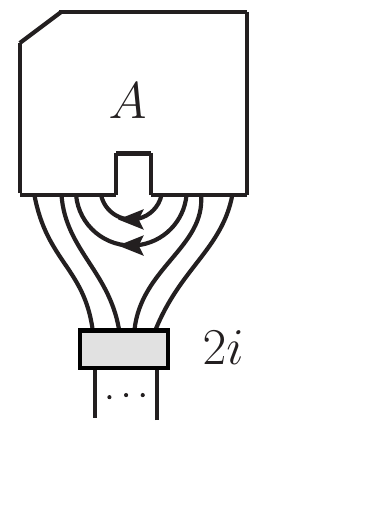}}
}
\ \ ,
\]
for all $0\le i\le n-1$. The resulting binary form is of degree $2i$ and the number of $\varepsilon$ arrows is $n-1-i$. 
we will also use the notation $\Pi_i(A)=A_i$ for components of this map.
In particular for $i=0$, and in view of (\ref{simtraceeq}), we see that the scalar $A_0$ is the trace of the endomorphism corresponding to $A$. Therefore, the Lie subalgebra $\mf{sl}_n$
is given by the kernel of the $\Pi_0$ projection map.

As an immediate consequence of the Clebsch-Gordan identities (\ref{CGSchureq}) and (\ref{CGdecompeq}),
we have that $\Pi$ is bijective, and its inverse is given by
\[
\Pi^{-1}(A_0,A_1,\ldots,A_{n-1})=
\sum_{i=0}^{n-1}
\frac{{\binom{n-1}{i}}^2}{\binom{n+i}{2i+1}} \times
\ \ \ 
\parbox{4.5cm}{\raisebox{-11ex}{
\includegraphics[width=4.5cm]{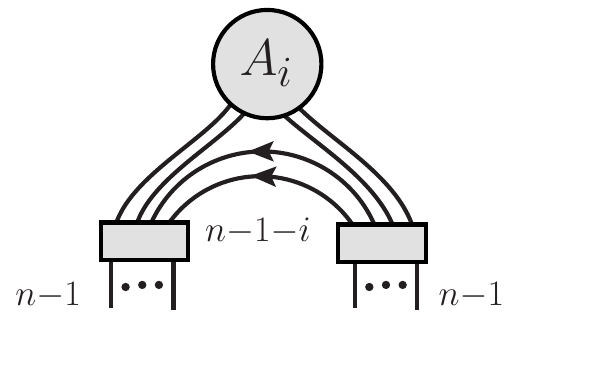}}
}
\ \ \ .
\]
It will be convenient to define, for $0\le j\le n-1$, the injective maps
\[
\mf{J}_i:\mc{H}_{2i}\longrightarrow \mc{H}_{0}\times\mc{H}_{2}\times\mc{H}_{4}\cdots\times\mc{H}_{2n-2}
\]
\[
A_i \longmapsto (0,\ldots,0,A_i,0,\ldots,0)\ .
\]

The main result of this section is the explicit computation of the composition $\circ$ when viewing elements of $\mf{gl}_n$, through the isomorphism $\Pi$, as sequences of binary forms $(A_0,A_1,\ldots,A_{n-1})$ of respective degree $0,2,4,\ldots,2n-2$.
For $0\le i,j,k\le n-1$, and binary forms $A_i\in\mc{H}_{2i}$ and $B_j\in\mc{H}_{2j}$,
 let
\[
\mathscr{C}_{k}^{i,j}(A_i,B_j):=
\Pi_k(\Pi^{-1}(\ \mf{J}_i(A_i))\circ\Pi^{-1}(\mf{J}_j(B_j))\ )\ ,
\]
which is a binary form of degree $2k$.
Since all our constructions are $\SL_2$-equivariant
and the Clebsch-Gordan decomposition for $\SL_2$ is multiplicity-free, it is clear that
the above expression is a numerical multiple of the transvectant $(A_i,B_j)_{i+j-k}$.
Our result gives an explicit formula for this coefficient, which essentially is a Wigner $6j$ symbol.

\begin{theorem}\label{prodthm}
We have
\[
\mathscr{C}_{k}^{i,j}(A_i,B_j)=
\bbone\{(i,j,k){\rm \ is \ a\ triad}\}
\times
\mathscr{P}_{k}^{i,j}
\times
(A_i,B_j)_{i+j-k}\ ,
\]
where
\begin{eqnarray}
\mathscr{P}_{k}^{i,j} & := & (-1)^{k+n-1}
\times\frac{i!\ j!\ k!}{(i+j-k)!(i+k-j)!(j+k-i)!}
\times\frac{\binom{n-1}{i}\binom{n-1}{j}}{\binom{n-1}{k}
\binom{n+i}{2i+1}\binom{n+j}{2j+1}} \nonumber\\
 & & \nonumber \\
 & & \times
\sum_q (-1)^{q}
\binom{q\! +\! 1}{i\! +\! j\! +\! k\! +\! 1}
\binom{i\! +\! j\! -\! k}{q\! -\! k\! -\! n\! +\! 1}
\binom{i\! +\! k\! -\! j}{q\! -\! j\! -\! n\! +\! 1}
\binom{j\! +\! k\! -\! i}{q\! -\! i\! -\! n\! +\! 1}\ . \nonumber \\
 & & \label{Pformulaeq}
\end{eqnarray}
\end{theorem}

In the above theorem, $\bbone\{\cdots\}$ denotes the indicator function of the condition between braces, the range for $q$ is $\bb{Z}$, with the usual convention that a binomial coefficient $\binom{s}{t}$ is zero unless $0\le t\le s$.

The proof of Theorem \ref{prodthm} is deferred to \S\ref{prfprodthmsec}.
From the theorem we immediately obtain similar formulas for the Lie bracket $[\cdot,\cdot]$ instead of the composition $\circ$, as follows.
We now work over the subalgebra $\mf{sl}_n$, and for $1\le i,j,k\le n-1$, as well as binary forms $A_i\in\mc{H}_{2i}$ and $B_j\in\mc{H}_{2j}$, we let
\[
\mathscr{D}_{k}^{i,j}(A_i,B_j):=
\Pi_k(\ [\ \Pi^{-1}(\mf{J}_i(A_i)),\Pi^{-1}(\mf{J}_j(B_j))\ ]\ )\ ,
\]
which is a binary form of degree $2k$.
We then obtain
\begin{equation}
\mathscr{D}_{k}^{i,j}(A_i,B_j)
=2\times \bbone\left\{\begin{array}{c}(i,j,k){\rm \ is \ a\ triad}\\
i+j-k{\rm \ is\ odd}
\end{array}\right\}
\times
\mathscr{P}_{k}^{i,j}
\times
(A_i,B_j)_{i+j-k}\ ,
\label{bracketeq}
\end{equation}
because of the symmetry $\mathscr{P}_{k}^{j,i}=\mathscr{P}_{k}^{i,j}$ and the identity $$(B_j,A_i)_{j+i-k}=(-1)^{i+j-k}(A_i,B_j)_{i+j-k}\ ,$$ which follows from the definition of transvectants. 

Let $\pi$ denote the the representation of $\SL_2$ on the space $\mc{H}_{n-1}$, i.e., $\pi(g)(F)=gF$, with the latter defined in (\ref{Ftransclassicaleq}), or (\ref{Ftransgraphicaleq}).
We will give a graphical formula for its derived representation $d\pi$ of $\mf{sl}_2$.
The latter will not be treated as the $n=2$ special case of the ongoing discussion with bihomogeneous forms, but simply as the space of $2\times 2$ matrices
\[
X=\begin{pmatrix}x_{11} & x_{12}\\ x_{21} & x_{22}\end{pmatrix}\ ,
\]
with $x_{11}+x_{22}=0$, and for which we will use the graphical notation
\[
\parbox{2.7cm}{\raisebox{-10ex}{
\includegraphics[width=2.7cm]{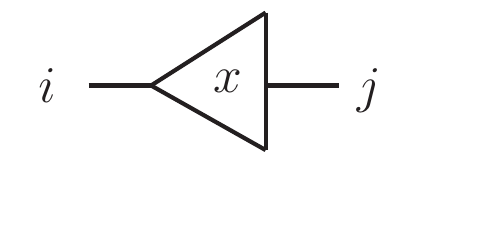}}
}
:= \ x_{ij}\ .
\]
We associate to $X$ a binary quadratic form $Q_X$ by
\[
\parbox{2.4cm}{\raisebox{-11ex}{
\includegraphics[width=2.4cm]{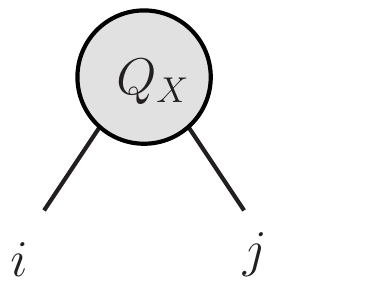}}
}
:= 
\parbox{3.4cm}{\raisebox{-10ex}{
\includegraphics[width=3.4cm]{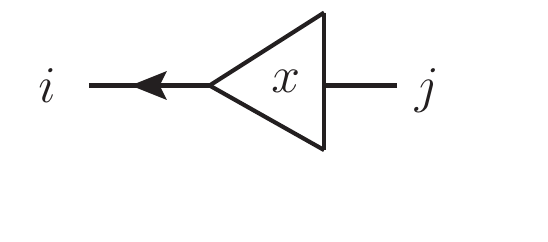}}
}
\ ,
\]
or, seeing the tensor of $Q_X$ as a symmetric $2\times 2$ matrix with the same name,
\[
Q_X:=\varepsilon X\ .
\]
The symmetry follows from the zero trace condition for $X$, since $(\varepsilon X)_{12}=x_{22}=-x_{11}=(\varepsilon X)_{21}$.

We will now use the definition $d\pi(X)=\left.\frac{d}{dt} \pi(e^{tX})\right|_{t=0}$, together with the identification of endomorphisms with bihomogeneous forms, in order to derive the following explicit graphical formula for $d\pi(X)$.

\begin{lemma}\label{Xpilem}
For all $X\in\mf{sl}_2$, we have
\begin{equation}
\parbox{2cm}{\raisebox{-13ex}{
\includegraphics[width=2cm]{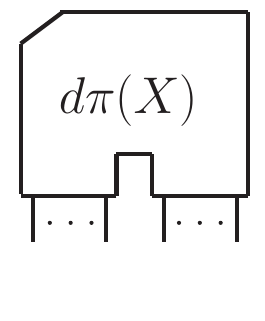}}
}
\ \ = (n-1)\times 
\ \ 
\parbox{5.8cm}{\raisebox{-11ex}{
\includegraphics[width=5.8cm]{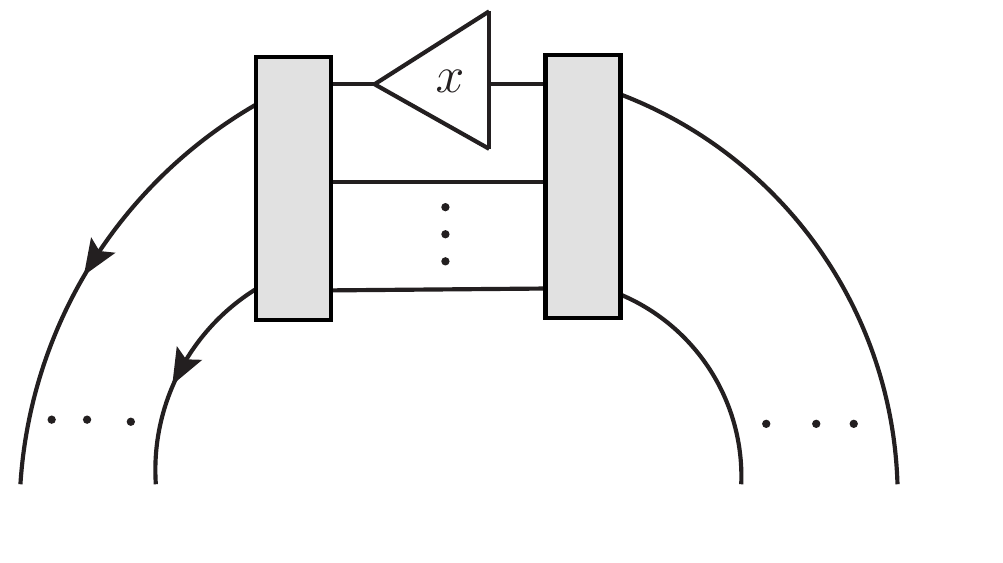}}
}
\ .
\label{graphicalpiXeq}
\end{equation}
Moreover, after processing through the map $\Pi$, we have for all $i$, $0\le i\le n-1$,
\[
\Pi_i(d\pi(X))=\delta_{i,1}\times \frac{n(n+1)}{6}\ Q_X\ .
\]
\end{lemma}

\noindent{\bf Proof:}
For $g=e^{tX}$, the trapeze-shaped tensor for the endomorphism $M=\pi(X)$ is
\[
\parbox{4cm}{\raisebox{-18ex}{
\includegraphics[width=4cm]{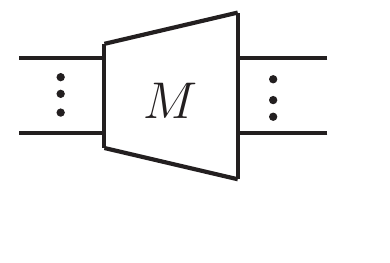}}
}
\ 
=\ \  
\parbox{4.5cm}{\raisebox{-12ex}{
\includegraphics[width=4.5cm]{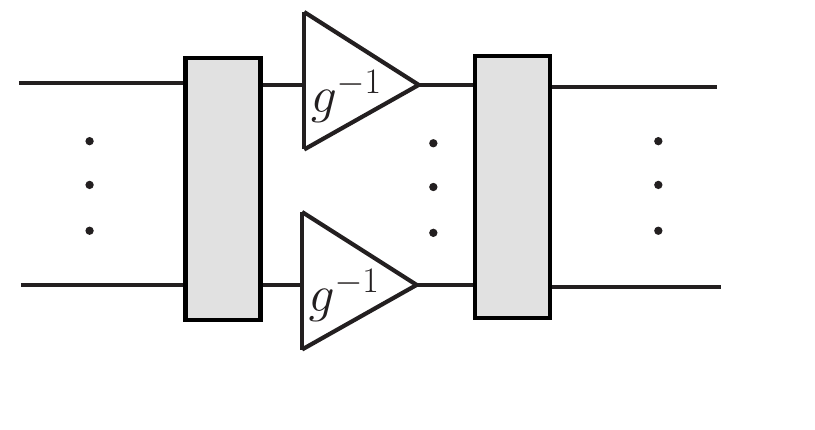}}
}
\ 
\ ,
\]
as results from (\ref{Ftransgraphicaleq}).
Note that having two symmetrizers is not necessary since one is enough, but this will be needed later for the linearization.
By (\ref{MtoAeq}), the corresponding bihomogeneous form is
\[ 
\parbox{2.4cm}{\raisebox{-18ex}{
\includegraphics[width=2.4cm]{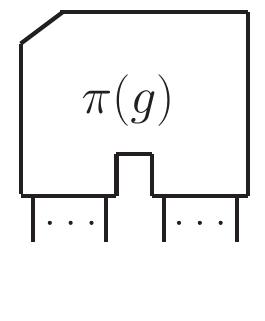}}
}
\ \ \ = \ \ 
\parbox{5.6cm}{\raisebox{-4ex}{
\includegraphics[width=5.6cm]{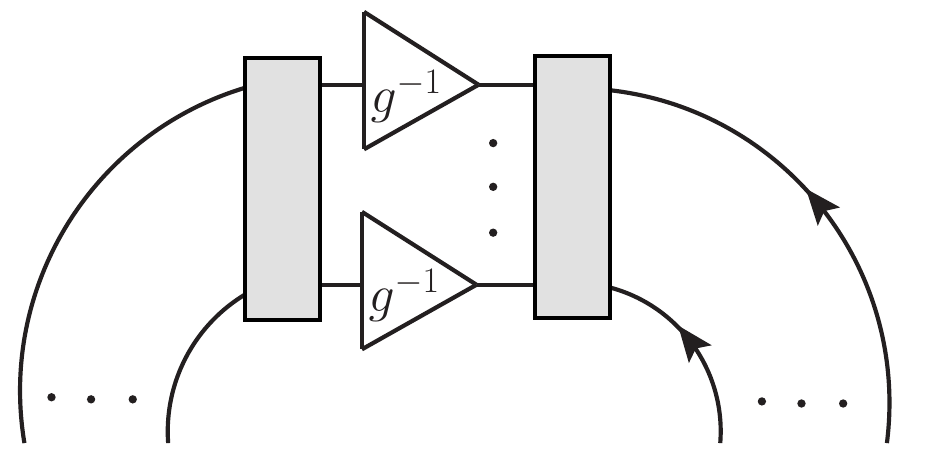}}
}
\]
\[
= \ \ 
\parbox{5.4cm}{\raisebox{-4ex}{
\includegraphics[width=5.4cm]{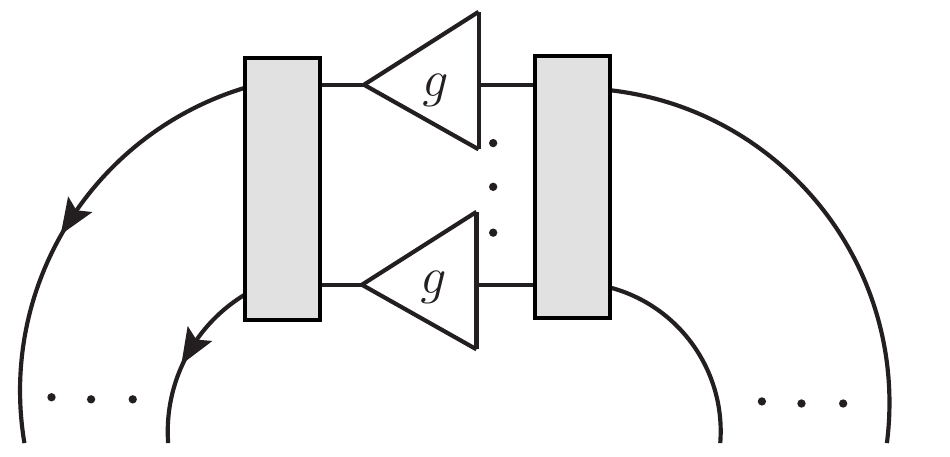}}
}
\ \ .
\]
The last equation follows from pushing the epsilon arrows through the right symmetrizer (see~\cite[Eq. (3.4)]{AbdesselamJKTR}), then using the graphical translation of the matrix identity $(g^{-1})^{\rm T}\varepsilon=\varepsilon g$ (Cramer's rule for $2\times 2$ matrices), and finally pushing the arrows through the left symmetrizer.
We now expand by multilinearity, and pick up the term linear in $X$ which is a sum of $(n-1)$ terms like the RHS of (\ref{graphicalpiXeq}), where the triangular piece for the matrix $X$ is placed on any one of the $n-1$ available horizontal strands between the symmetrizers.
The symmetrizers allow the exchange of these ``ladder rungs'', and therefore (\ref{graphicalpiXeq}) holds.
As a result, concerning the second part of the lemma, we have
\[
\Pi_i(d\pi(X))=(n-1)\times \ \ 
\parbox{4.6cm}{\raisebox{-4ex}{
\includegraphics[width=4.6cm]{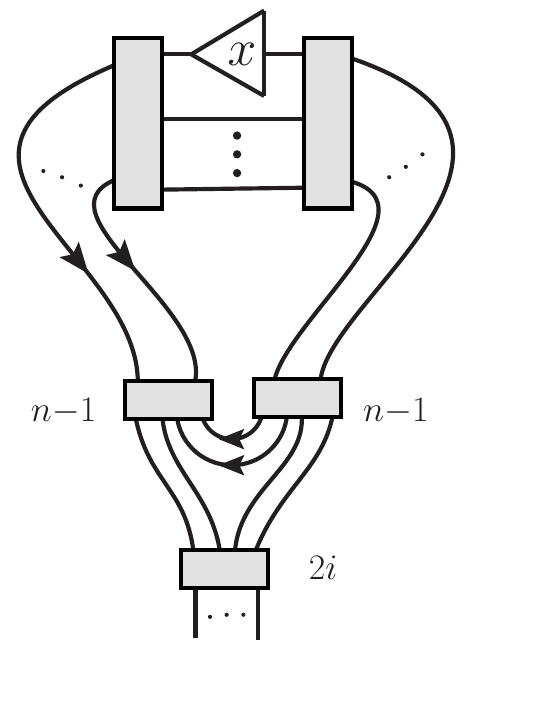}}
}
\]
\[
=(n-1)\times \ \ 
\parbox{4.6cm}{\raisebox{-4ex}{
\includegraphics[width=4.6cm]{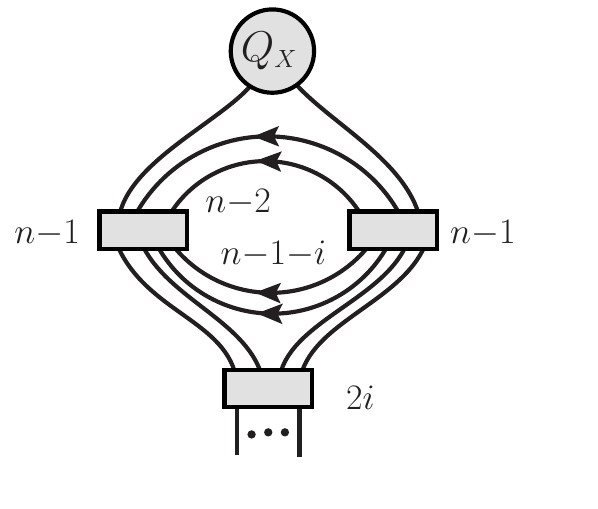}}
}
\]
\[
=\delta_{i,1}\times\frac{n(n+1)}{6} \times \ \ 
\parbox{2.4cm}{\raisebox{-11ex}{
\includegraphics[width=2.4cm]{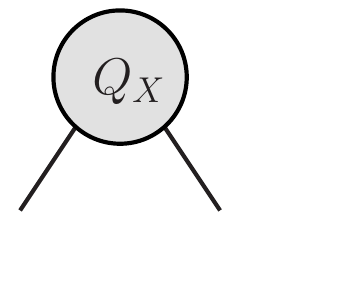}}
}
\ \ .
\]
The first line is the definition of the projection $\Pi_i$ together with (\ref{graphicalpiXeq}) we just proved.
In passing from the first to the second line, we used the idempotence of symmetrizers to remove the top two, and we also used the definition of $Q_X$.
The last line results from the graphical Schur Lemma (\ref{CGSchureq}) which applies, because $Q_X$ is symmetric, i.e., one can insert a size two symmetrizer on top, at no cost.
\qed

We can now give a concrete incarnation for our main protagonists $E,H,F\in\mf{gl}_n$
seen as bihomogenous forms.
Starting from the $\mf{sl}_2$ generators
\[
e=\begin{pmatrix}
0 & 1\\
0 & 0
\end{pmatrix}
\ \ ,\ \ 
h=\begin{pmatrix}
1 & 0\\
0 & -1
\end{pmatrix}
\ \ ,\ \ 
f=\begin{pmatrix}
0 & 0\\
1 & 0
\end{pmatrix}\ ,
\]
we just let
\[
E:=d\pi(e)\ \ ,\ \ 
H:=d\pi(h)\ \ ,\ \ 
F:=d\pi(f)\ ,
\]
with corresponding quadratic binary forms $Q_X$ therefore given by
\[
Q_e(\mathbf{x})=-x_{2}^{2}\ \ ,\ \ 
Q_h(\mathbf{x})=-2x_{1}x_{2}\ \ ,\ \ 
Q_f(\mathbf{x})=x_{1}^{2}\ .
\]
Of course, $E,H,F$ form a standard $\mf{sl}_2$-triple.
Indeed, $d\pi$ is a Lie algebra morphism. As an opportunity to practice the graphical calculus,
we invite the reader
to prove this property diagrammatically, starting with the pictures in (\ref{graphicalpiXeq}) as definitions.
The basic idea, is that when taking the commutator in $\mf{gl}_n$ of $d\pi(X)$ and $d\pi(Y)$,
one gets pictures which contain a portion such as
\[
\parbox{3.8cm}{\raisebox{-11ex}{
\includegraphics[width=3.8cm]{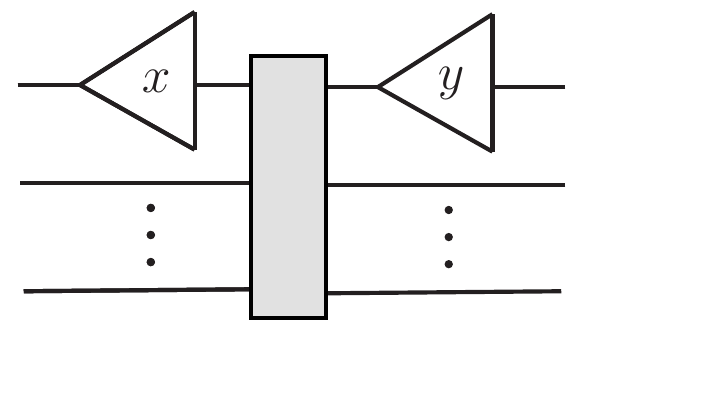}}
}
\ \ .
\]
When expanding the middle symmetrizer as a sum of permutations, with probability $\frac{1}{n-1}$, the $x$ and $y$ triangles get mounted in series,
and with remaining probability $\frac{n-2}{n-1}$, they get mounted in parallel. However, the parallel placement cancels when combining both terms of the commutator. The series placement, on the other hand rebuilds a small triangle for the $2\times 2$ matrix $[X,Y]$.
The same reasoning allows one to easily compute powers of $E$ and $F$, where the reverse phenomenon happens: only the parallel placement survives because the matrices $e,f$ are nilpotent. 

We now define, for $0\le i\le n-1$, the subspace $V_i=\Pi^{-1}(\mf{J}_i(\mc{H}_{2i}))$ of $\mf{gl}_n$. We clearly have $\mf{gl}_n=\bigoplus_{0\le i\le n-1}V_i$,
and  $\mf{sl}_n=\bigoplus_{1\le i\le n-1}V_i$. One expects these decomposition to coincide with the decomposition into irreducible modules for the principal $\mf{sl}_2$ considered earlier in Section \ref{Sec:natural-basis}, however, this requires a few routine checks which one could do using graphical calculations as suggested in the previous paragraph.
We will instead use the following proposition which is a specialization of (\ref{bracketeq})
and which will be useful in the proof of Proposition \ref{Gcompprop}.

\begin{proposition}\label{adcompprop}
For $1\le i\le n-1$, $0\le j\le n-1$, $A\in V_1$ and $B\in V_i$, we have
\[
\Pi_j(\ad_A(B))=-\delta_{ij}\times\frac{12\ i}{n(n+1)}\times
(\Pi_1(A),\Pi_i(B))_1\ .
\]
\end{proposition}

\noindent{\bf Proof:}
We apply (\ref{bracketeq}), and note that the triad condition gives $|i-1|\le j\le i+1$
which, together with the parity condition that $1+i-j$ is odd, implies that the result is zero unless $j=i$. Then, we have
$\Pi_i(\ad_A(B))=2\times \mathscr{P}_{i}^{1,i}\times 
(\Pi_1(A),\Pi_i(B))_1$, so the result reduces to the computation of $\mathscr{P}_{i}^{1,i}$.
The sum in (\ref{Pformulaeq}) now only contains two terms, namely, the one with $q=n+i-1$ and the one with $q=n+i$.
More precisely, we get
\[
\mathscr{P}_{i}^{1,i}=\frac{i!^2}{(2i-1)!}\times\frac{6}{n(n+1)}\times\frac{1}{\binom{n+i}{2i+1}}
\left[\binom{n+i}{2i+2}\binom{2i-1}{i-1}
-\binom{n+i+1}{2i+2}\binom{2i-1}{i}\right]\  .
\]
Using 
\[
\binom{2i-1}{i-1}=\binom{2i-1}{i}\ \ ,\ \ 
\binom{n+i+1}{2i+2}=\binom{n+i}{2i+2}+\binom{n+i}{2i+1}\ , 
\] 
and cleaning up, we obtain $\mathscr{P}_{i}^{1,i}=-\frac{6i}{n(n+1)}$, and we are done.
\qed

The proposition immediately implies, as expected, that $\mf{sl}_n$, seen as a module for the imbedded $\mf{sl}_2$ Lie subalgebra generated by the triple $(E,H,F)$ defined above, decomposes into the submodules $V_1,\ldots,V_{n-1}$.

We now give explicit formulas for the powers of $E$ and $F$, and for the $G_{i,j}$ basis defined in Equation \eqref{def-G}.

\begin{lemma}\label{Epowerlem}
For all $k\ge 0$, and all $i$ with $0\le i\le n-1$, we have
\[
\Pi_i(E^k)=\delta_{i,k}\ (-1)^k\ k!\times\frac{\binom{n+k}{2k+1}}{\binom{n-1}{k}}\times x_2^{2k}\ .
\]
\end{lemma}

\noindent{\bf Proof:}
We proceed by induction on $k$. For $k=0$, $E^k=\mathrm{id}$ or rather the associated bihomogeneous form graphically given by
\[
\parbox{4.8cm}{\raisebox{-11ex}{
\includegraphics[width=4.8cm]{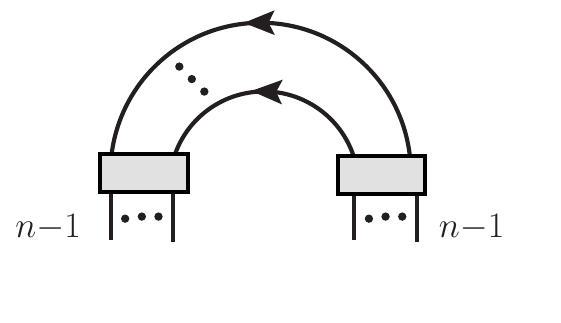}}
}
\ \ .
\]
The special case $k=\ell=0$ of
the identity (\ref{CGSchureq}) then shows the LHS, a scalar, is equal to $\delta_{i,0}\times n$,
where the $n$ factor accounts for the trace of the identity map on $\mc{H}_{n-1}$.
The RHS agrees, which establishes the base case of the induction.
For the $k=1$ case needed for the induction step, we note that since $Q_e=-x_2^2$,
Lemma \ref{Xpilem} immediately gives
\[
\Pi_i(E)=
-\delta_{i,1}\times\frac{n(n+1)}{6}\times x_2^2\ ,
\]
which agrees with the RHS after simplification and confirms that $E\in V_1$.
Now assume the result is true for $k$. Then, $E\in V_k$ and
Theorem \ref{prodthm}
shows that
$\Pi_i(E\circ E^k)$ is zero unless the triad condition $|k-1|\le i\le k$ holds, in which case $\Pi_i(E\circ E^k)$ is proportional to the transvectant $(x_2^2, x_2^{2k})_{1+k-i}$. From formula (\ref{transveq}), one immediately sees that this is zero unless the transvectant is of zero-th order, i.e., $i=k+1$. Hence, $E^{k+1}\in V_{k+1}$ if $k\le n-2$ and is zero if $k>n-2$. For $k\le n-2$, we only need to compute $\Pi_{k+1}(E\circ E^k)$, in order to complete the proof by induction.
By the induction hypothesis, the $k=1$ case, and Theorem \ref{prodthm}, we have
\[
\Pi_{k+1}(E\circ E^k)=
(-1)\times\frac{n(n+1)}{6}
\times
(-1)^k\ k!\times\frac{\binom{n+k}{2k+1}}{\binom{n-1}{k}}
\times \mathscr{P}_{k+1}^{1,k}
\times (x_2^2,x_2^{2k})_0\ .
\]
The zero-th transvectant is just the ordinary product, and the sum (\ref{Pformulaeq}) for $\mathscr{P}_{k+1}^{1,k}$ reduces to the single term with $q=k+n$. After simplification, we thus obtain the RHS of the formula in the lemma, with $k+1$ instead of $k$.
\qed

\begin{lemma}\label{Fpowerlem}
For all $k\ge 0$, and all $i$ with $0\le i\le n-1$, we have
\[
\Pi_i(F^k)=\delta_{i,k}\ k!\times\frac{\binom{n+k}{2k+1}}{\binom{n-1}{k}}\times x_1^{2k}\ .
\]
\end{lemma}

The proof is entirely similar to that of the previous lemma and is thus left to the reader.

We now explicitly compute, in classical invariant-theoretic fashion, the basis $G_{i,j}$, labeled by the integers $i,j$ subject to the ranges $1\le i\le n-1$, $-i\le j\le i$, for the Lie algebra $\mf{sl}_n$ seen inside $\mc{H}_{n-1}\otimes\mc{H}_{n-1}$. Recall that, by definition, $G_{i,j}= (\ad_F)^{i-j}(E^i)$.

\begin{proposition}\label{Gcompprop}
The element $G_{i,j}$ belongs to $V_i$ and its only nonzero $\Pi$ projection is
\[
\Pi_i(G_{i,j})=(-1)^j\times i!\times
\frac{\binom{n+i}{2i+1}}{\binom{n-1}{i}}
\times\frac{(2i)!}{(i+j)!}\times x_1^{i-j}x_2^{i+j}\ .
\] 
\end{proposition}

\noindent{\bf Proof:}
We first consider the specialization of Proposition \ref{adcompprop} to $A=F$. For  $B\in V_i$, and denoting $\Pi_i(B)$ by $B_i$, the proposition and the $k=1$ case of Lemma \ref{Fpowerlem} give us
\begin{eqnarray*}
\Pi_i(\ad_F(B)) & = & -\frac{12i}{n(n+1)}
\times
\left(\frac{n(n+1)}{6}x_1^2,B_i\right)_1\\
 & = & -2i\times\ \ 
\parbox{4cm}{\raisebox{-13ex}{
\includegraphics[width=4cm]{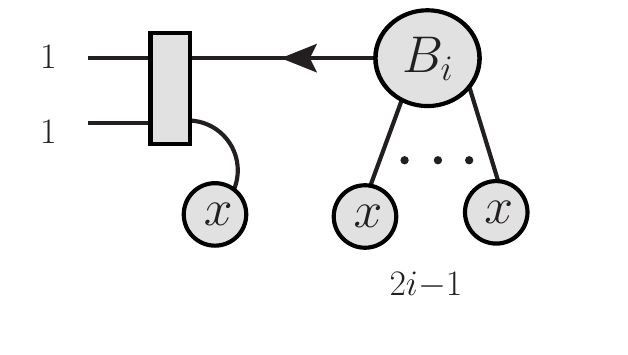}}
}
 \\
 & = & -2i\times x_1\times\  
\parbox{2.9cm}{\raisebox{-15ex}{
\includegraphics[width=2.9cm]{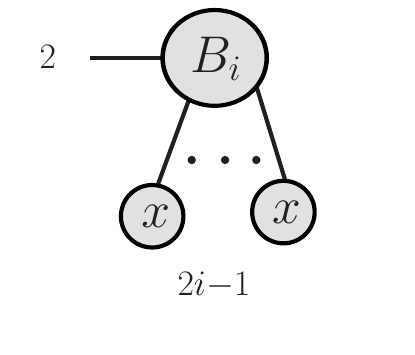}}
}
\\
 & = & -x_1\frac{\partial B_i}{\partial x_2}(\mathbf{x})\ .
\end{eqnarray*}
In the first picture above, the two loose legs on the left carry indices set equal to the value 1. In the second picture, the single loose leg on the left of the $B_i$ ``blob'' carries an index set equal to 2. Of course, one could also compute the transvectant non-graphically using formula (\ref{transveq}), in order to arrive at the same result.
We now have, thanks to Lemma \ref{Fpowerlem},
\begin{eqnarray*}
\Pi_i(G_{i,j}) & = & \left(-x_1\frac{\partial}{\partial x_2}\right)^{i-j}
\ \left[
(-1)^i\ i!\times\frac{\binom{n+i}{2i+1}}{\binom{n-1}{i}}\times x_2^{2i}
\right] \\
 & = & (-1)^j\ i!\times\frac{\binom{n+i}{2i+1}}{\binom{n-1}{i}}\times x_1^{i-j} 
\left(\frac{\partial}{\partial x_2} \right)^{i-j}
x_2^{2i}\ ,
\end{eqnarray*}
since the operators of multiplication by $x_1$ and derivation with respect to $x_2$ commute. Finishing the computation of the derivative, and cleaning up, we recover the desired formula.
\qed

Since the monomials in $x_1,x_2$ give a basis for $\mc{H}_{2i}$, and the last proposition shows that the $G_{i,j}$, produced by repeated application of the lowering operator $\ad_F$, are essentially proportional to these monomials, we completed the checks needed to verify that the $\mf{sl}_2$-modules earlier defined by $V_i=\Pi^{-1}(\mf{J}_i(\mc{H}_{2i}))$ are indeed irreducible. Since there are $n-1$ of them in the decomposition $\mf{sl}_n=\bigoplus_{1\le i\le n-1}V_i$, i.e., as much as the rank of $\mf{sl}_n$, we also checked that the triple $(E,H,F)$ produced in this section is a principal $\mf{sl}_2$ triple for $\mf{sl}_n$. 

\subsection{The proof of Theorem \ref{prodthm}}
\label{prfprodthmsec}

Referring to the settings and notations from Section \ref{Sec:graphical-calc-for-sln}, let
\[
C_k:=\mathscr{C}_{k}^{i,j}(A_i,B_j):=
\Pi_k(\ \Pi^{-1}(\mf{J}_i(A_i))\circ\Pi^{-1}(\mf{J}_j(B_j))\ )\ .
\]
The ``blob'' of $C_k$ is given by
\[
\parbox{2cm}{\raisebox{-13ex}{
\includegraphics[width=2cm]{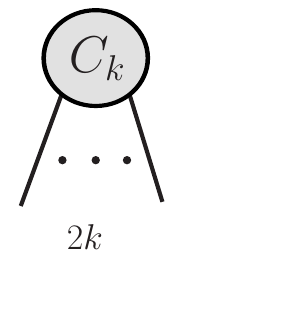}}
}
\ \ = 
\frac{{\binom{n-1}{i}}^2{\binom{n-1}{j}}^2}{\binom{n+i}{2i+1} \binom{n+j}{2j+1}}
\times\ \ 
\parbox{8cm}{\raisebox{-11ex}{
\includegraphics[width=8cm]{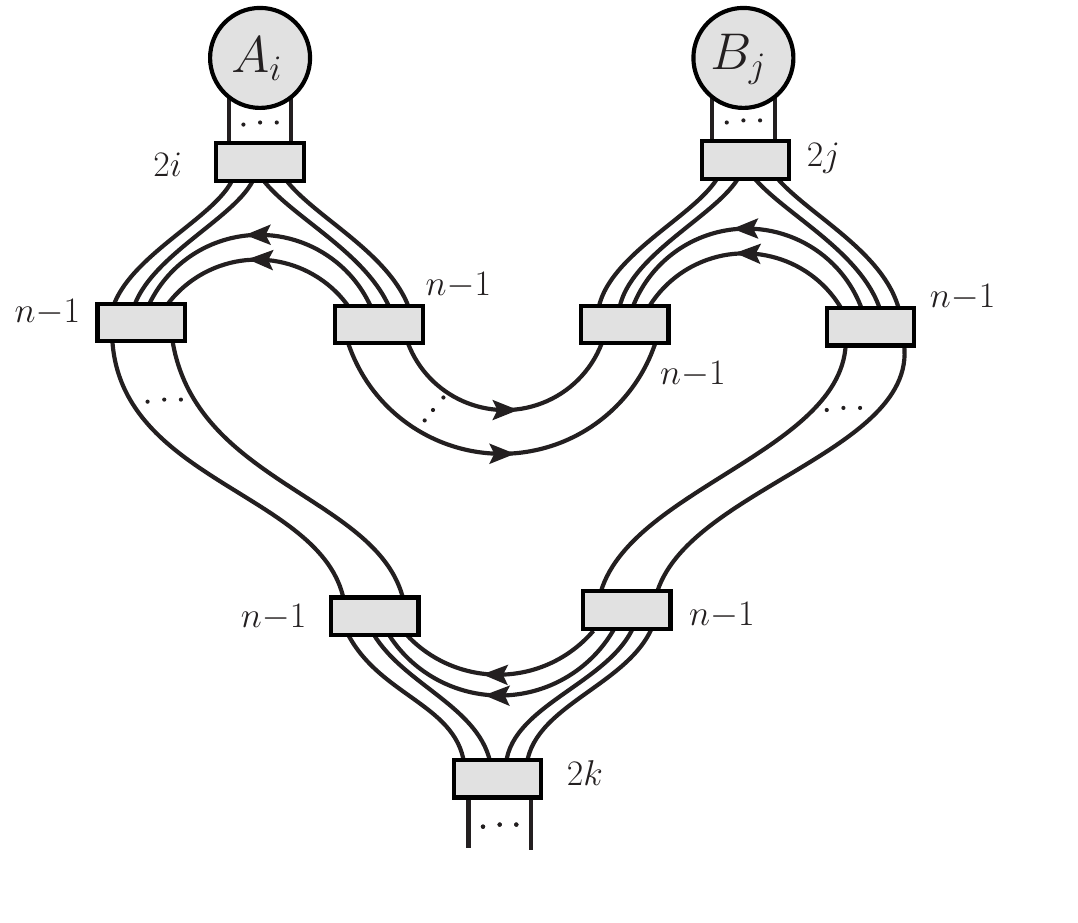}}
}
\ \ ,
\]
where the composition of maps can be read by ``scanning through'' the picture from top to bottom. We put near each symmtrizer an indication of its size. We also put redundant symmetrizers of respective sizes $2i$ and $2j$ immediately under the ``blobs'' for $A_i$ and $B_j$. This is because it is where we will insert the Clebsch-Gordan identity (\ref{CGdecompeq}), which results in the equation
\[
\parbox{1.5cm}{\raisebox{-13ex}{
\includegraphics[width=1.5cm]{pic47.pdf}}
}
\!\!\!\! =\frac{{\binom{n-1}{i}}^2{\binom{n-1}{j}}^2}{\binom{n+i}{2i+1} \binom{n+j}{2j+1}}
\times\sum_{\ell=0}^{\min(2i,2j)}
\frac{\binom{2i}{\ell}\binom{2j}{\ell}}{\binom{2i+2j-\ell+1}{\ell}}
\times \!\!
\parbox{7.3cm}{\raisebox{-11ex}{
\includegraphics[width=7.3cm]{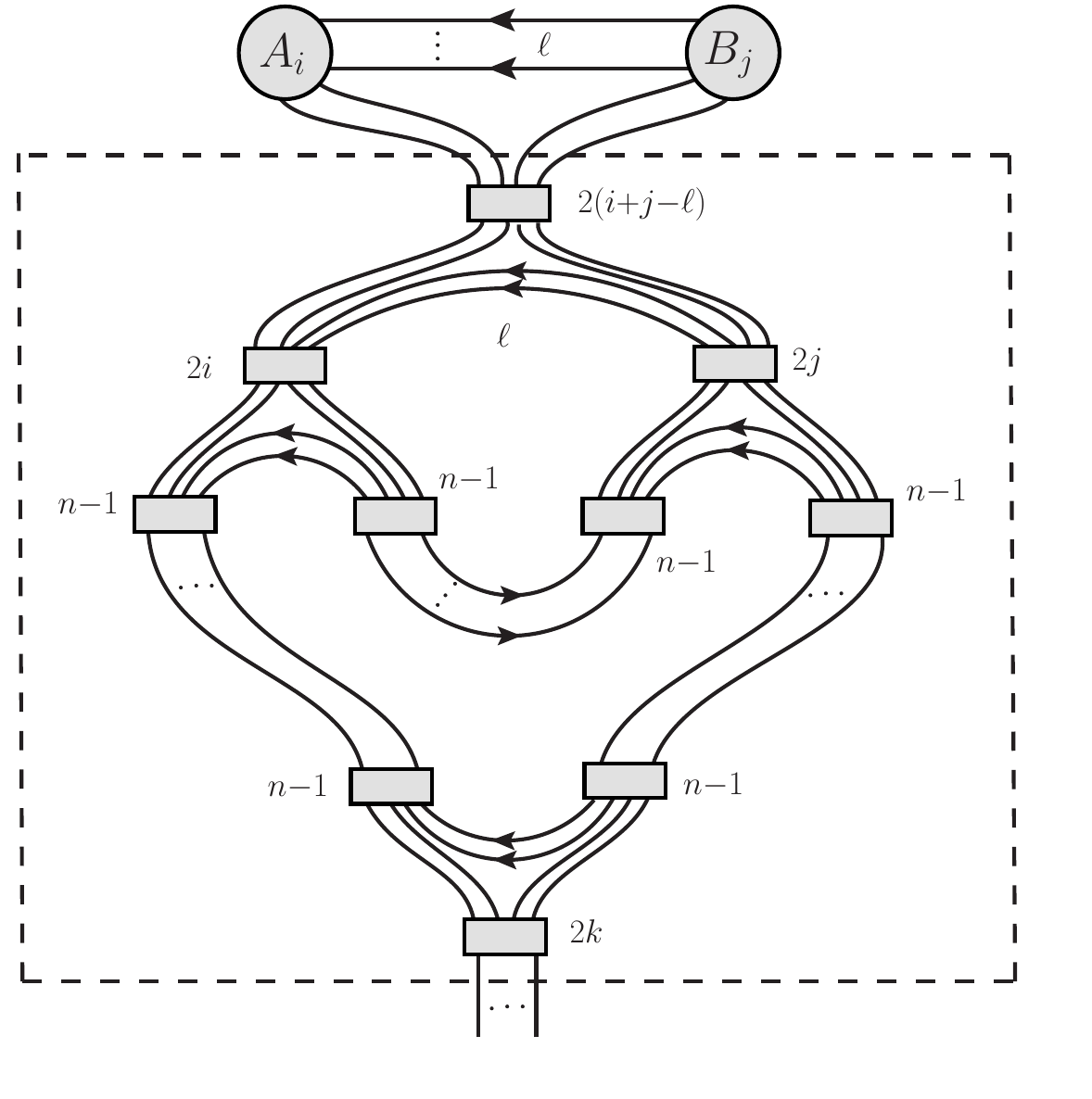}}
}
\ \ .
\]
The portion inside the dotted box, read from top to bottom is an $\SL_2$-equivariant map $\mc{H}_{2(i+j-\ell)}\rightarrow\mc{H}_{2k}$. The graphical Schur Lemma~\cite[Prop. 3.2]{AbdesselamJKTR}
forces $2(i+j-\ell)=2k$, i.e., the contributing value of the summation index is $\ell=i+j-k$. The allowed range for $\ell$ gives the indicator function of $(i,j,k)$ being a triad.
For $\ell=i+j-k$, Schur's Lemma also gives the proportionality
\[
\Gamma:= 
\parbox{8.5cm}{\raisebox{-11ex}{
\includegraphics[width=8.5cm]{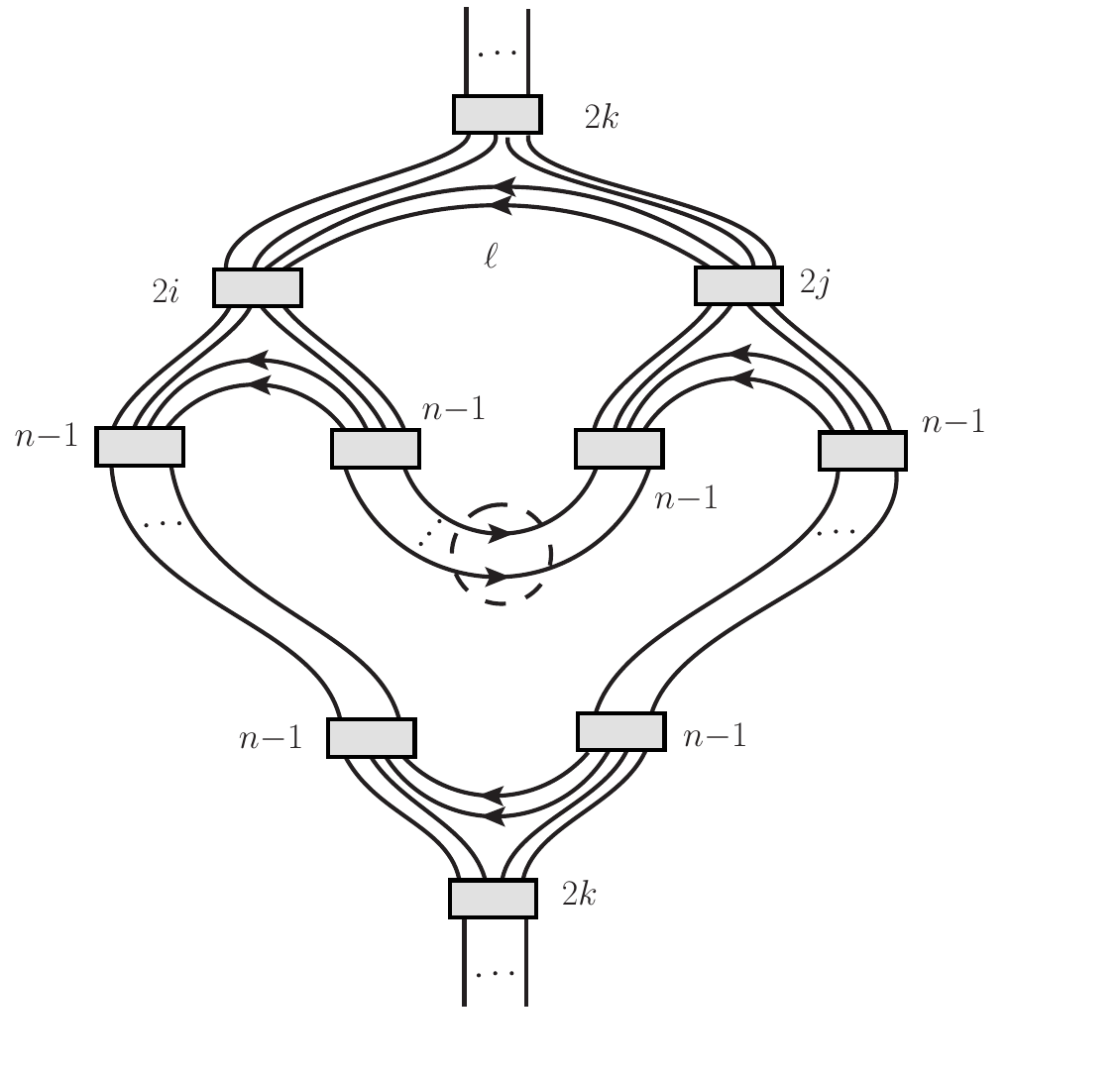}}
} =\gamma \times 
\parbox{2cm}{\raisebox{-8ex}{
\includegraphics[width=2cm]{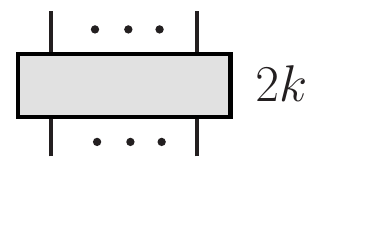}}
}
\ ,
\]
for a suitable of the factor $\gamma$.
We thus have
\begin{equation}
\mathscr{C}_{k}^{i,j}(A_i,B_j) =\bbone\{(i,j,k){\rm \ is \ a\ triad}\}
\times
\frac{{\binom{n-1}{i}}^2{\binom{n-1}{j}}^2}{\binom{n+i}{2i+1} \binom{n+j}{2j+1}}
\times
\frac{\binom{2i}{i+j-k}\binom{2j}{i+j-k}}{\binom{i+j+k+1}{i+j-k}}
\times
\gamma
\times (A_i,B_j)_{i+j-k}\ .
\label{Cintermeq}
\end{equation}
We then determine $\gamma$ by relating the graphical expression $\Gamma$ to (\ref{6jdefeq}).
For this, we push the $\varepsilon$ arrows singled out by the dotted circle through the symmetrizer on the left, as in~\cite[Eq. (3.4)]{AbdesselamJKTR}.  Out of these $n-1$ arrows, the $n-1-i$ leftmost ones will disappear because of the $2\times 2$ matrix identity $\varepsilon \varepsilon^{\rm T}=I$. We then flip the direction of the remaining $i$ rightmost ones so they point towards the symmetrizer of size $2i$ above them. This produces a $(-1)^i$ factor, 
namely, we get
\[
\Gamma =(-1)^i \times\ \ \ \ 
\parbox{8.5cm}{\raisebox{-11ex}{
\includegraphics[width=8.5cm]{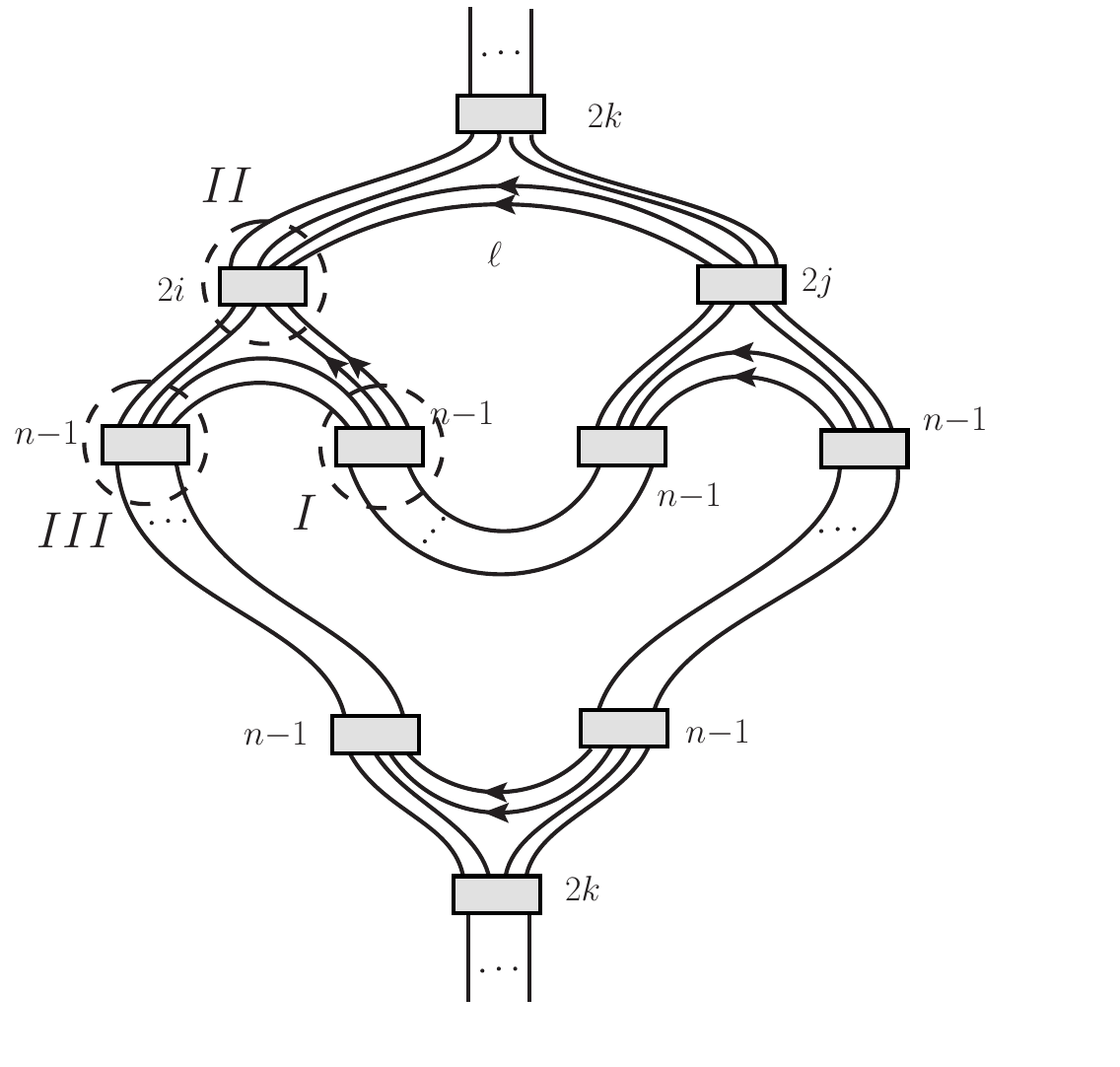}}
}
\ .
\]
The next three moves I, II, III, use the idempotence
of symmetrizers in both directions. We let the symmetrizer in the dotted circle I be absorbed by the symmetrizer to its right.
We do the reverse operation and duplicate the symmetrizer II. Finally we let the symmetrizer III be absorbed by the one below it.
This gives
\[
\Gamma=(-1)^i\times
\parbox{8.5cm}{\raisebox{-11ex}{
\includegraphics[width=8.5cm]{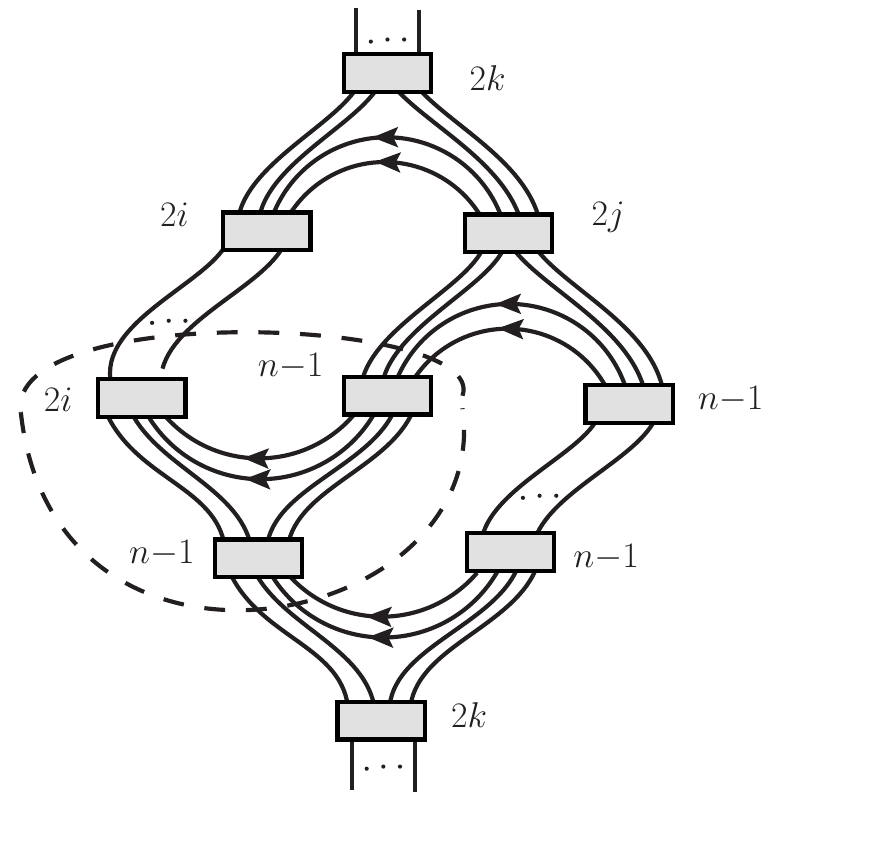}}
}
\ .
\]
Essentially, the previous three symmetrizers in the dotted circles have been rotated 
(around the pivot symmetrizer I), and produced the triangular structure inside the new dotted curve.
Comparing with (\ref{6jdefeq}), we read off
\[
\gamma = (-1)^i\times
\rho\left[
\begin{array}{ccc}
i & \frac{n-1}{2} & \frac{n-1}{2}\\
\frac{n-1}{2} & k & j
\end{array}
\right]\ .
\]
We use these entries in order to invoke (\ref{rhoeq}), and multiply the result by the other factors from (\ref{Cintermeq}). After a rather tedious simplification and reorganization of factorials into binomial coefficients, we arrive at the formula (\ref{Pformulaeq}) for the quantity $\mathscr{P}_{k}^{i,j}$.

\begin{Remark}
The coefficients $\mathscr{P}_{k}^{i,j}$ experimentally tend to be nonzero, in general,
except for $\mathscr{P}_{2}^{2,2}$ with $n=4$. Their vanishing, or not, is equivalent to that of the $6j$ symbol
\[
\left\{\begin{array}{ccc}
i & \frac{n-1}{2} & \frac{n-1}{2}\\
\frac{n-1}{2} & k & j
\end{array}\right\}
=
\left\{\begin{array}{ccc}
i & j & k \\
\frac{n-1}{2} & \frac{n-1}{2} & \frac{n-1}{2}\\
\end{array}\right\}\ ,
\]
using the previously mentioned symmetries to write it in a more memorable form.
See the last section of~\cite{AbdesselamRMS}, for a conjecture about the nonvanishing of a subset of these $6j$ symbols.
\end{Remark}

\subsection{Structure constant computations}

We are ready to prove the main result of this article, namely Theorem \ref{Thm:2} which we give again here.

\begin{theorem}\label{explicitomegathm}
For the Lie algebra $\mathfrak{sl}_n$, the structure constants $\omega_{k,\ell}^{(m)}$ vanish for $m$ even, and  are given, when $m$ is odd, by
\[
\omega_{k,\ell}^{(m)}=\mathcal{Q}_{k,\ell}^{(m)}\times\mathcal{R}_{k,\ell}^{(m)}\ ,
\] 
with
\[
\mathcal{Q}_{k,\ell}^{(m)}:=2\times(-1)^{k+\ell+n-1}\times(2k+2\ell-2m+1)\times
\frac{k!^2\ \ell!^2\ (n-k-\ell+m-1)!}{m!\ (n+k+\ell-m)!}\ ,
\]
and
\[
\mathcal{R}_{k,\ell}^{(m)}:=
\sum_{q\in\mathbb{Z}}(-1)^q
\binom{q\! +\! 1}{2k\! +\! 2\ell \! -\! m\! +\! 1}\!
\binom{m}{q\! -\! k\! -\! \ell\! +\! m\! -\! n\! +\! 1}\!
\binom{2k\! -\! m}{q\! -\! \ell\! -\! n\! + \! 1}\!
\binom{2\ell\! -\! m}{q\! -\! k\! -\! n\! +\! 1}\ ,
\]
with the convention that binomials are defined as zero if lying outside of Pascal's triangle.
\end{theorem}

To prove the theorem, we use the following equation (see Equation \eqref{Eq:useful-eq}) with $1\le k,\ell\le n-1$, which can be taken as definition of the structure constants:
\[
[E^k,F^{\ell}]=
\sum_{m=1}^{2\min(k,\ell)}
\frac{(-1)^{\ell}\omega_{k,\ell}^{(m)}}{(2\ell-m)!}\ 
G_{k+\ell-m,k-\ell}\ .
\]

\noindent{\bf Proof:}
We apply (\ref{bracketeq}), for $1\le i\le n-1$, in order to compute the projection $\Pi_i([E^k,F^{\ell}])$,
and see that the result is zero unless $(k,\ell,i)$ form a triad and $k+\ell-i$ is odd. This gives the range of the new index $m:=k+\ell-i$, namely, $0\le m\le 2\min(k,\ell)$ and 
shows that $\omega_{k,\ell}^{(m)}=0$ unless
$m$ is odd. 
We assume these conditions hold, in what follows.
Using (\ref{bracketeq}), with $i=k+\ell-m$, the formulas in Lemmas \ref{Epowerlem} and \ref{Fpowerlem},
as well as the easy transvectant computation
\[
(x_2^{2k},x_1^{2\ell})_m=(-1)^m x_1^{2\ell-m} x_2^{2k-m}\ ,
\]
we obtain
\[
\Pi_{k+\ell-m}([E^k,F^{\ell}])=\alpha\times x_1^{2\ell-m} x_2^{2k-m}\ ,
\]
with
\begin{equation}
\alpha:=2\times \mathscr{P}_{k+\ell-m}^{k,\ell}\times (-1)^k \ k!
\ \frac{\binom{n+k}{2k+1}}{\binom{n-1}{k}}
\times \ell!\ \frac{\binom{n+\ell}{2\ell+1}}{\binom{n-1}{\ell}}
\times (-1)^m\ .
\label{alphaintermeq}
\end{equation}
On the other hand, by Proposition \ref{Gcompprop}, we have
\[
\Pi_{k+\ell-m}(G_{k+\ell-m,k-\ell})=\beta \times x_1^{2\ell-m} x_2^{2k-m}\ ,
\]
with
\begin{equation}
\beta:=(-1)^{k-\ell}
\ (k+\ell-m)!\times
\frac{\binom{n+k+\ell-m}{2k+2\ell-2m+1}}{\binom{n-1}{k+\ell-m}}
\times\frac{(2k+2\ell-2m)!}{(2k-m)!}\ .
\label{betaintermeq}
\end{equation}
Comparing coefficients of the $G$ basis vectors, we
obtain
\[
\omega_{k,\ell}^{(m)}=(-1)^{\ell} \ (2\ell-m)!\times\frac{\alpha}{\beta}\ .
\]
We substitute (\ref{alphaintermeq}) and (\ref{betaintermeq}), and insert the formula for $\mathscr{P}_{k+\ell-m}^{k,\ell}$
from (\ref{Pformulaeq}), and after considerable simplification, we arrive at the desired result.
\qed

Note that the literature on Wigner's $3nj$ symbols, and $6j$ symbols in particular, is quite extensive. Some computer algebra systems also have preprogrammed functions for the evaluation of $6j$ symbols with the standard definition recalled in (\ref{6jformulaeq}).
For instance, in Mathematica, the relevant command which performs this evaluation is ``SixJSymbol''.
Therefore, for the reader's convenience, we record below the formula for the $\omega$'s directly in terms of standard $6j$ symbols:
\begin{eqnarray}\label{Eq:6j-expr}
\omega_{k,\ell}^{(m)} &=& (-1)^{k+\ell+n-1}\times 2\times(2k+2\ell-2m+1)\times\frac{k!\ \ell!}{(k+\ell-m)!} \nonumber\\
 & & \times\sqrt{\frac{(2k-m)!(2\ell-m)!(n+k)!(n+\ell)!(n-k-\ell+m-1)!}{m!(n-k-1)!(n-\ell-1)!(2k+2\ell-m+1)!(n+k+\ell-m)!}} \nonumber\\
 & & \times\left\{
 \begin{array}{ccc}
 k & \ell & k+\ell-m\\
 \frac{n-1}{2} & \frac{n-1}{2} & \frac{n-1}{2}
 \end{array}
 \right\}
\end{eqnarray}
The latter follows from comparing (\ref{6jformulaeq}) and the formula in Theorem \ref{explicitomegathm}.

\section{Other simple complex Lie algebras}\label{Sec:other-types}

Some aspects from the previous sections carry over to a general simple complex Lie algebra $\g$. Fix a principal $\mathfrak{sl}_2$-triple in $\g$. The decomposition of $\g$ into irreducible $\mathfrak{sl}_2$-modules is given by $\g\cong \oplus_i V_i$ with $\dim V_i = 2m_i+1$. The $m_i$ are known as the \emph{exponents} of $\g$. 
We gather some insights about the structure constants by distinguishing the type of $\g$.

\subsection{Type B and C}
For $\g=\mathfrak{so}_{2n+1}$ or $\g=\mathfrak{sp}_{2n}$, there is a natural inclusion into type $A$ given by $\mathfrak{so}_{2n+1}\subset \mathfrak{sl}_{2n+1}$ and $\mathfrak{sp}_{2n}\subset \mathfrak{sl}_{2n}$. 

We can characterize the elements of $\mathfrak{so}_{2n+1}$ or $\mathfrak{sp}_{2n}$ as the fixed point set of two Lie algebra involutions.
For type $B$, the involution is $i_B(M)=-M^T$ and for type $C$ it is $i_C(M)=-JM^tJ^{-1}$ where $J$ is the matrix of the symplectic form.

The special property in type $B$ and $C$ is that a principal $\mathfrak{sl}_2$-triple $(E,H,F)$ is also principal in the ambient Lie algebra of type $A$. Therefore, we obtain the theory in type $B$ (resp. $C$) via the $i_B$-invariant (resp. $i_C$-invariant) part of the theory of type $A$.

In both cases, the decomposition $\g\cong \bigoplus_k V_k$ runs now over $k$ odd. The highest weight vectors can be chosen to be the odd powers of $E$. 
Hence for $k, \ell, m$ odd we have
\begin{equation}\label{Eq:type-B}
    \omega_{k,\ell}^{(m)}(\mathfrak{so}_{2n+1}) = \omega_{k,\ell}^{(m)}(\mathfrak{sl}_{2n+1})
\end{equation}
and the same for type $C$. For $k$ or $\ell$ even, the structure constant is zero.

\subsection{Type D}

For $\g=\mathfrak{so}_{2n}$, there is also an inclusion into $\mathfrak{sl}_{2n}$. But a principal $\mathfrak{sl}_2$-triple of $\mathfrak{so}_{2n}$ is not principal in the ambient space $\mathfrak{sl}_{2n}$.

The decomposition into irreducible $\mathfrak{sl}_2$-representations is given by
$$\mathfrak{so}_{2n}\cong V_{n-1}\oplus \bigoplus_{k=1}^{n}V_{2k-1}.$$
We have odd indices (as for type $B$) and one additional index which we denote by $n'$ ($n'=n-1$, but for $n$ even it should not be mixed with the odd index $n-1$).

Fix the following principal nilpotent element (all empty entries are zero):
\begin{equation}
F=\left(\begin{array}{@{}cccc|cccc@{}}
0&&&&&&& \\
1&0&&&&&& \\
&\ddots&\ddots&&&&& \\
&&1&0&&&& \\ \hline
&&1&0&0&&& \\
&&&-1&-1&0&& \\
&&&&&\ddots &\ddots& \\
&&&&&&-1&0\\
\end{array}\right).
\end{equation}

The centralizer $Z(F)$ is generated by odd powers of $F$ and a special matrix $S$ given by $S=E_{n,1}-E_{n+1,1}+E_{2n,n}-E_{2n,n+1}$:

\begin{equation}\label{matrixS}
S=\left(\begin{array}{@{}ccc|ccc@{}}
  &&&&& \\
	&&&&&\\
	1&&&&& \\ \hline
	-1&&&&&\\
	&&&&&\\
	&&1&-1&&
\end{array}\right)
\end{equation}

The matrix $S$ has no apparent link to $F$ or $E$, so it seems unnatural to attribute a monomial to it (in particular for $n$ even in which case we have twice the summand $V_{n-1}$). Two observations can be made: 
\begin{proposition}
    The following  structure constants vanish in type $D_n$ for all $k,\ell,m$: 
    \begin{equation}\label{Eq:type-D}
    \omega_{k,\ell}^{(n')}=0, \; \omega_{k,n'}^{(m)}=0 \;\text{ and }\; \omega_{n',n'}^{(n')}=0.
    \end{equation}
    In addition, for the ``usual'' odd indices $k,\ell,m$, we have 
    \begin{equation}\label{Eq:type-D-2}
    \omega_{k,\ell}^{(m)}(\mathfrak{so}_{2n})=\omega_{k,\ell}^{(m)}(\mathfrak{so}_{2n-1})=\omega_{k,\ell}^{(m)}(\mathfrak{sl}_{2n-1}).
    \end{equation}
\end{proposition}

The idea of the proof uses the folding procedure from type $D_n$ to type $B_{n-1}$. 

\begin{proof}
For the first part, consider the Lie algebra involution $i$ on $\mathfrak{so}_{2n}$ corresponding to the horizontal symmetry in the Dynkin diagram. On the level of the simple roots $(\alpha_1,..., \alpha_{n-1},\alpha_n)$, it acts via $i(\alpha_k)=\alpha_k$ for $1\leq k\leq n-2$, $i(\alpha_{n-1})=\alpha_n$ and $i(\alpha_n)=\alpha_{n-1}$.
It is easy to check that $i(F)=F$ and $i(S)=-S$. The fixed point set of $i$ is generated by $\bigoplus V_k$ and the anti-invariant part by $V_{n'}$.
By the involution property, we immediately get Equation \eqref{Eq:type-D}.

The fixed point set of $i$ is isomorphic to the Lie algebra $\mathfrak{so}_{2n-1}$. Since the odd powers of $F$ are in the fixed point set, we get $\omega_{k,\ell}^{(m)}(\mathfrak{so}_{2n})=\omega_{k,\ell}^{(m)}(\mathfrak{so}_{2n-1})$ and we conclude with Equation \eqref{Eq:type-B}.
\end{proof}

\subsection{Exceptional types}

For the remaining five exceptional types, it would be possible to compute all the structure constants, once one has fixed the highest weight vectors. We have only done the computation in the smallest case, the type $G_2$. 

For type $G_2$, the decomposition of $\g_2$ into irreducible $\mathfrak{sl}_2$-representations is $\g_2 \cong V_1\oplus V_5$. In the representation of lowest dimension, which is 7, the highest weight vectors can be chosen to be $E$ and $E^5$. So a natural choice is to attribute the monomial $x^{10}/10!$ to $E^5$. A direct computation (using the nice article \cite{wildberger}) gives the two non-trivial structure constants:
$$\omega_{5,5}^{(9)}(\mathfrak{g}_2) = -\frac{27000}{7} \;\;\text{ and }\;\; \omega_{5,5}^{(5)}(\mathfrak{g}_2) = -\frac{30}{7}.$$

\appendix

\section{Tables with structure constants}\label{Sec:appendix}

We list here all structure constants for $\mathfrak{sl}_n$ with $n\leq 6$. From the definition, we know that $\omega_{k,\ell}^{(m)}=\omega_{\ell,k}^{(m)}$, so we consider only $k\geq \ell$. Also it directly follows from the definition that $\omega_{k,1}^{(m)}=\delta_{m,1}$. Hence we only consider $k\geq \ell\geq 2$.

For $\mathfrak{sl}_3$, the only non-trivial structure constant is $\omega_{2,2}^{(3)}= -2$.
	
The structure constants for $\mathfrak{sl}_4$ are given by:

\vspace{0.3cm}
\begin{center}
\renewcommand{\arraystretch}{1.3}
\begin{tabular}{|l|l|}
\hline 
$\omega_{2,2}^{(3)} = -24/5$ & $\omega_{2,2}^{(1)}=2/5$ \\
\hline
$\omega_{3,2}^{(3)}=-3$ & 0\\
\hline
$\omega_{3,3}^{(5)} = 54/5$ &  $\omega_{3,3}^{(3)} = 3/5$\\
\hline
\end{tabular} 
\end{center}
\vspace{0.3cm}

The structure constants for $\mathfrak{sl}_5$  are given by:

\vspace{0.3cm}
\begin{center}
\renewcommand{\arraystretch}{1.3}
\begin{tabular}{|l|l|}
\hline 
$\omega_{2,2}^{(3)} = -42/5$ & $\omega_{2,2}^{(1)}=2/5$ \\
\hline
$\omega_{3,2}^{(3)}=-48/7$ & $\omega_{3,2}^{(1)}=3/14$\\
\hline
$\omega_{3,3}^{(5)} = 216/5$ &  $\omega_{3,3}^{(3)} = -6/5$\\
\hline
$\omega_{4,2}^{(3)} = -4$ &  0\\
\hline
$\omega_{4,3}^{(5)} = 144/7$ &  $\omega_{4,3}^{(3)} = 6/7$\\
\hline
$\omega_{4,4}^{(7)} = -576/5$ &  $\omega_{4,4}^{(5)} = -24/5$\\
\hline
\end{tabular} 
\end{center}
\vspace{0.3cm}

The structure constants for $\mathfrak{sl}_6$  are given by:

\vspace{0.3cm}
\begin{center}
\renewcommand{\arraystretch}{1.3}
\begin{tabular}{|l|l|l|}
\hline 
$\omega_{2,2}^{(3)} = -64/5$ & $\omega_{2,2}^{(1)}=2/5$ & 0\\
\hline
$\omega_{3,2}^{(3)}=-81/7$ & $\omega_{3,2}^{(1)}=3/14$ &0\\
\hline
$\omega_{3,3}^{(5)} = 3888/35$ &  $\omega_{3,3}^{(3)} = -17/5$ & $\omega_{3,3}^{(1)} = 1/14$\\
\hline
$\omega_{4,2}^{(3)} = -80/9$ &  $\omega_{4,2}^{(1)} = 2/15$ & 0\\
\hline
$\omega_{4,3}^{(5)} = 540/7$ &  $\omega_{4,3}^{(3)} = -3/7$ & 0\\
\hline
$\omega_{4,4}^{(7)} = -4608/7$ &  $\omega_{4,4}^{(5)} = 16/3$ & $\omega_{4,4}^{(3)} = 16/21$\\
\hline
$\omega_{5,2}^{(3)} = -5$ & 0  & 0\\
\hline
$\omega_{5,3}^{(5)} = 100/3$ &  $\omega_{5,3}^{(3)} = 1$ & 0\\
\hline
$\omega_{5,4}^{(7)} = -1800/7$ & $\omega_{5,4}^{(5)} = -60/7$  & 0\\
\hline
$\omega_{5,5}^{(9)} = 14400/7$ & $\omega_{5,5}^{(7)} = 200/3$  & $\omega_{5,5}^{(5)} = 20/21$\\
\hline
\end{tabular} 
\end{center}
\vspace{0.3cm}

\bibliographystyle{plain}
\bibliography{ref}

\end{document}